\documentclass[12pt,reqno]{amsart}
\input{epsf.sty}
\catcode `\@=11
\def\numberbysection{\@addtoreset{equation}{section}
         \renewcommand{\theequation}{\thesection.\arabic{equation}}}
\numberbysection
\def\subsubsection{\@startsection{subsubsection}{3}%
  \normalparindent{.5\linespacing\@plus.7\linespacing}{-.5em}%
  {\normalfont\bfseries}}
\def\cal{\mathcal}

\def\del{\partial}

\def\b{\beta}
\def\g{\gamma}

\def\eps{\epsilon}
\def\l{\lambda}

\def\Cal{\cal}

\def\Z2{{\bf Z}/2{\bf Z}}
\def\Zn{{\bf Z}/n{\bf Z}}
\def\Zn+1{{\bf Z}/(n+1){\bf Z}}

\def\la{\langle}
\def\ra{\rangle}
\def\Z2{{\bf Z}/2{\bf Z}}
\def\Zn{{\bf Z}/n{\bf Z}}
\def\Znn{{\bf Z}/(n+1){\bf Z}}

\begin{document}

\title{Orbifolding Frobenius Algebras}

\author[Ralph M. Kaufmann]{Ralph M. Kaufmann\\
University of Southern California, Los Angeles, USA$^*$}
\email{kaufmann@math.usc.edu}
\thanks{$^*$ Partially supported by NSF grant
DMS \#0070681}
\address{University of Southern California, Los Angeles, USA}

\begin{abstract}
We study the general theory of Frobenius algebras with group actions.
These structures arise when one is studying the algebraic structures associated to
a geometry stemming from a physical  theory with a global finite gauge group, i.e.\ orbifold theories.
In this context, we introduce and axiomatize these algebras.
Furthermore, we define geometric cobordism categories whose functors to the category of
 vector spaces are parameterized
by these algebras. The theory is also extended to the graded and
super--graded cases. As an application, we consider Frobenius
algebras having some additional properties making them more tractable. These properties
are present in Frobenius algebras arising as quotients of Jacobian ideal, such as those having their
origin in quasi--homogeneous singularities and their symmetries.
\end{abstract}

\maketitle

\section*{Introduction}

The subject of this exposition is the general theory of Frobenius algebras with group actions.
These structures arise when one is studying the algebraic structures
 stemming from a geometry associated to a physical
theory with a global finite gauge group [DW,DVVV,IV,V].
The most prominent example of this type in mathematics is the Gromov--Witten theory of
orbifolds [CR], which are global quotients. The use of orbifold constructions is the cornerstone
of  the original mirror construction [GP]. The orbifolds under study in that context are so--called
Landau--Ginzburg orbifold theories, which have so far  not been studied
mathematically. These correspond to the Frobenius manifolds coming from singularities and
are studied as examples in detail in the present paper.

A common aspect of the physical treatment of quotients by group actions is the appearance of
so called twisted sectors. This roughly means that if one wishes to take the quantum version
of quotient by a group action, one first has to construct an object for each element in the group together
with a group action on this object and in a second step take invariants in all of these components.
Whereas classically one considers only the $G$--invariants of the original object which constitutes the
sector associated to the identity, the untwisted sector.

We give the first complete axiomatic treatment, together with a natural geometric interpretation for
this procedure, which provides a common basis for applications such as singularities with
symmetries or Laundau--Ginzburg orbifolds, orbifold cohomology and quantum cohomology of global
quotients and in a sense all other string--orbifold versions of classical theories.

Our treatment shows that the construction of twisted sectors  is not  merely an
auxiliary artifact, but is essential. This is clearly visible in the case of the  Frobenius algebra
associated to the singularity of type $A_n$ together with a $\Znn$ action and the singularity
of type $A_{2n-3}$ with $\Z2$ action, which are worked out in
detail in the last paragraph. In particular the former example exhibits a version of mirror symmetry
in which  it is self--dual. The twisted sectors  are the key in this mirror duality, since it is
the sum of twisted sectors that is dual to the untwisted one.

In the present work, we develop the theory of orbifold Frobenius algebras along the now classical
lines of Atiyah, Dubrovin, Dijkgraaf and Segal.
 I.e.\ we start by introducing the algebraic structures in an
axiomatic fashion.
There is an important difference to the theory with trivial group in that there
are two  structures with slightly different $G$-action to be considered which are nevertheless present.
These versions differ by a twist with a character.
In the singularity version, on the non--twisted level,
they correspond to the different $G$--modules
on the cohomology and
the Milnor ring [Wa]. In physical terms, these two structures are related by spectral flow.
To be even more precise, one structure carries the natural multiplication and the other
the natural scalar product.
The next step is a cobordism realization of the theory. There are ramifications of the present exposition
with [FQ, T], if the above mentioned character is trivial. This is, however, not necessarily the case
and in the most interesting examples worked out in the last paragraph, this is not the case.

The key structure here is a non--commutative multiplication on the sum of all the twisted sectors,
even before taking invariants. This is a novel approach never
considered before.

In order to apply our theory to Landau--Ginzburg models or singularity
theory --in a sense the original building blocks for mirror symmetry--,
we introduce a large class of examples, so--called special $G$--Frobenius algebras.

This class contains the class of
Jacobian $G$--Frobenius algebras which in turn encompass the
singularity examples and those of manifolds whose cohomology ring can be described as a quotient by a
Jacobian ideal. Here it is important to note that everything can be done in a super (i.e.\ $\Z2$--graded)
 version. This introduces a new degree of freedom into the construction corresponding to the choice
of parity for the twisted sectors.
Lastly, we explicitly work out several examples including the transition from
the singularity $A_{2n+3}$ to $D_n$ via a quotient by $\Z2$, this is
the first purely mathematical version of this correspondence avoiding
path integrals.
In this situation the presence of the Ramond algebra explains the
results obtained by [Wa] in the situation of singularities with
symmetries as studied by Arnold.
We furthermore show that $*/G$ leads to (twisted)-group algebras.

For Jacobian Frobenius algebras, we also introduce a duality
transformation which allows us to show that orbifolding  plays
the role of mirror symmetry. We show that
the pair $(A_{n},A_{1})$ is mirror dual to $(A_{1},A_{n})$ via
orbifolding by $\Znn$. Applying the same principle to $A_{1}/\Znn$ we
find the underlying Frobenius algebra structure for the A--model
realization of $A_{n}$ [W,JKV].

In the case that the Frobenius algebra one starts out with comes
from a semi--simple Frobenius manifold
and the quotient of the twisted sector is not trivial, there
is  unique extension to the level
of Frobenius manifolds. This is the case in the above example of  $A_{2n+3}$ and $D_n$.

The general theory presented here applies to the orbifold cohomology
of global quotients, as worked out in [FG], who found our postulated
non--commutative structure in that situation.

Furthermore the theory of special G--twisted Frobenius algebras
in the case where G is the symmetric group sheds new light on the construction
of [LS], explaining the uniqueness and how the general structure
of their multiplication is fixed [K3]. Recently the structure
of discrete torsion was uncovered and explained in [K4] to be given by
forming tensor products in the sense introduced below with
twisted group rings. The synthesis of these results and an
application to the Hilbert scheme can be found in [K5].

\section*{Acknowledgements}
This work is the written version of talk presented at the
WAGP2000 conference in October 2000,
which contained all algebraic constructions and examples in full detail,
and the talk at the workshop on ``Mathematical aspects of
orbifold string theory'' in May 2001, which contained the
geometric cobordism part of the third paragraph.
I wish to thank the organizers of these conferences, 
which were immensely stimulating.

I would also like to thank the Max--Planck Institut f\"ur
Mathematik in Bonn where much of this work
was completed during my stays in 2000 and 2001 as well as
the IHES in Bures--sur--Yvette.  The research has profited from conversations
with Yu.\ Manin, Y.\ Ruan, V.\ Turaev, C.\ Hertling and B.\ Penner all of
whom it  is a pleasure to thank.

Finally, I would like to acknowledge the support from the NSF grant 
DMS\#0070681.

\section{Frobenius Algebras and Cobordisms}
In this paragraph, we recall the  definition of a Frobenius algebra and its relation to
the cobordism--category definition of a topological field theory [A,Du,Dij].

\subsection {Frobenius algebras}
\subsubsection {Definition}  A {\em Frobenius algebra} (FA) over
a field  $k$ of characteristic 0 is
$\la A,\circ, \eta, 1\ra$, where

\begin{tabular}{ll}
A&is a finite dim $k$--vector space, \\
$\circ$&is a multiplication on $A$:  $\circ: A\otimes A\rightarrow A$,\\
$\eta$&is a non-degenerate bilinear form on $A$, and\\
$1$&is a fixed element in $A$ -- the unit.
\end{tabular}

 satisfying the following axioms:

\begin{itemize}

\item[a)] {\em Associativity}

$(a\circ b) \circ c =a\circ (b \circ c)$
\item[b)] {\em Commutativity}

$a\circ b =b\circ a$
\item[c)]
{\em Unit}:

$\forall \in a: 1 \circ a =  a \circ 1 = a$

\item[d)]
{\em Invariance}:

$\eta(a,b\circ c) = \eta(a\circ b,c)$

\end{itemize}
\subsubsection{Remark}
\label{convert}
By using $\eta$ to identify $A$ and $A^*$ -- the dual vector space of $A$--
these objects define a one--form $\eps \in A^*$ called the co--unit  and
 a three--tensor $\mu \in A^*\otimes A^*\otimes A^*$.

Using dualization and invariance these data are interchangeable
with $\eta$ and $\circ$ via the following formulas. Explicitly,
after fixing a basis $(\Delta_{i})_{i \in I}$ of $A$, setting $\eta_{ij}:= \eta(\Delta_i, \Delta_j)$ and denoting the
inverse metric by $\eta^{ij}$\\

$\eps(a):=\eta(a,1)$,

$\mu(a,b,c) :=\eta(a\circ b,c)=\eps(a\circ b\circ c)$,

$a \circ b = \sum_{ij}\mu(a,b,\Delta_{i})\eta^{ij}\Delta_{j}$ and

$\eta(a,b)=\eps(a\circ b)$.

\subsubsection{Notation} We call $\rho \in A$ the element dual to $\eps$. This is the element which is
Poincar\'e dual to $1$.

\subsubsection{Grading}

A graded Frobenius algebra is a Frobenius algebra together with a
group grading of the vector space $A$: $A=\bigoplus_{i \in I} A_i$
where $I$ is a group together with the following compatibility
equations: denote the $I$--degree of an element by $\deg$.
\begin{itemize}

\item[1)] $1$ is homogeneous; $1 \in A_d$ for some $d\in I$
\item[2)] $\eta$ is homogeneous of degree $d+D$.\\
 I.e.\ for homogeneous elements $a,b$
$\eta (a,b)=0$ unless $\deg(a)+\deg(b)= d+D$.

This means that $\eps$ and $\rho$ are of degree $D$.

\item[3)] $\circ$ is of degree $d$, this means that $\mu$ is of degree $2d+D$,
where again this means that

$\deg(a\circ b) = \deg(a)+\deg(b) -d$
\end{itemize}

\subsubsection{Definition} An even derivation $E\in Der(A,A)$ of a $G$ twisted
Frobenius algebra $A$ is called an {\em Euler field}, if it is conformal and
is natural w.r.t. the multiplication, i.e. for some $d,D\in k$ it satisfies:

\begin{equation}
\eta(Ea,b)+\eta(a,Eb)= D\eta(a,b)
\end{equation}
and
\begin{equation}
E(ab) = Ea \, b +a \, Eb - d \, ab
\end{equation}

Such derivation will define a grading on $A$  by its set of eigenvalues.

\subsubsection{Remark}
For this type of grading, we will use the group $\bf Q$.
There are two more versions of grading: 1) a ${\bf Z}/2{\bf Z}$ super-grading, which will be discussed
in \ref{supergrading} and 2) a grading by  a finite group $G$, which is the content of section \ref{orb}.

\subsubsection{Definition}
Given an I--graded Frobenius algebra $A$, we define its {\em  characteristic series} as

\begin{equation}
\chi_{A}(t) := \sum_{i \in I} \dim (A_i) t^i
\end{equation}

We refer to the set $\{d,D;i:\dim A_i \neq 0\}$ as the (I--)spectrum  of $A$.
\subsubsection{Scaling}
If the group indexing the grading has the structure of a $\Lambda$-module, where $\Lambda$ is a ring,
 we can scale the grading by an element $\lambda \in \Lambda$.
We denote the scaled Frobenius algebra by
$\lambda A:= \bigoplus_{i\in I}A_{\lambda i}$ where $(\lambda A)_i = A_{\lambda i}$

\begin{equation}
\chi_{\lambda A}(t) := \sum_{i \in I} \dim (\lambda A)_i t^i =
t^{\lambda}\chi_A(t)
\end{equation}

 It is sometimes --but not always-- convenient to normalize in such a way that $\deg(1)=0$ where
 $0$ is the unit in $I$.
In the case that the grading is given by $\bf Q$, this means
$\deg_{\bf Q}(1)=0$ and in the finite group case $\deg_G(1)=e$
where now $e$ is the unit element of $G$.

\subsubsection{Operations}
There are two natural operations on Frobenius algebras, the direct sum and the tensor product.
Both of these operation extend to the level of Frobenius manifolds,
while the generalization of the
direct sum is straightforward, the generalization of the tensor product to the level of Frobenius
manifolds is quite
intricate [K1].

\medskip
\noindent Consider two Frobenius algebras ${\mathcal A}'=  \la A', \circ',\eta', 1' \ra$ and
${\mathcal A}'' = \la A'', \circ'',$ $ \eta'', 1'' \ra$ .

\subsubsection{Direct sum}

Set ${\mathcal A}'\oplus {\mathcal A}'' :=  \la A,\circ,\eta, 1
\ra =
  \la A'\oplus A'', \circ' \oplus \circ'',\eta' \oplus \eta'',1' \oplus 1'' \ra$.
E.g. $ (a',a'') \circ (b',b'') = (a'\circ' b', a''\circ'' b'')$.
The unit is $1=1'\oplus 1''$ and the co-unit is $\eps= (\eps',\eps'')$.

\subsubsection{Lemma} If both Frobenius algebras are graded by the same $I$ then their direct sum
inherits a natural grading if and only if the gradings can be scaled s.t.
\begin{equation}
\label{scale}
D'+d'=D''+d'':= D+d
\end{equation}
where
\begin{equation}
\label{scalecond}
D=D'=D''
\end{equation}
in {\em this} scaling.

Furthermore, the unit will have degree $d'=d''=d$.

{\bf Proof}

The equation  (\ref{scale}) ensures
that the three tensor $\mu$ is homogeneous of degree $D+2d$.
The homogeneity of $\eta$ yields the second condition:
for $\eta$ to be homogeneous, we must have that after scaling
$\rho'$ and $\rho''$ are homogeneous of the same degree $D'=D''=D$.
The two equations together imply the homogeneity of $1=(1',1'')$ of degree $d=d'=d''$.

\subsubsection{Tensor product}

Set ${\mathcal A}'\otimes {\mathcal A}'' :=  \la A,\circ,\eta, 1 \ra =
  \la A'\otimes A'', \circ' \otimes \circ'',\eta' \otimes \eta'',1' \otimes 1'' \ra$.\\
E.g. $ (a',a'') \circ (b',b'') = (a'\circ' b', a''\circ'' b'')$.
The unit is $1=1'\otimes 1''$ and the co-unit is  $\eps= \eps'\otimes \eps''$.

There are {\em no conditions} for grading. I.e.\ if both Frobenius algebras are $I$--graded
there is a natural induced $I$--grading on their tensor product. The unit is of degree $d=d'+d''$,
the co--unit has degree $D'+D''$ and the multiplication is homogeneous of degree $d=d'+d''$.

\subsubsection{Super-grading}
\label{supergrading}
For an element $a$ of a super vector space $A= A_0 \oplus A_1$ denote by $\tilde a$ its $\Z2$ degree,
i.e.\ $\tilde a = 0$ if it is even ($a\in A_0$) and $\tilde a = 1$ if it is odd ($a\in A_1$).

\subsubsection{Definition}
A {\em super Frobenius algebra} over
a field  $k$ of characteristic 0 is
$\la A, \circ, \eta, 1\ra$, where

\begin{tabular}{ll}
A&is a finite dim $k$--super vector space, \\
$\circ$&is a multiplication on $A$:  $\circ: A\otimes A\rightarrow A$, \\
&which preserves the $\Z2$--grading\\
$\eta$&is a non-degenerate even bilinear form on $A$, and\\
$1$&is a fixed even element in $A_0$ -- the unit.
\end{tabular}

 satisfying the following axioms:

\begin{itemize}

\item[a)] {\em Associativity}

$(a\circ b) \circ c =a\circ (b \circ c)$
\item[b)] {\em Super--commutativity}

$a\circ b = (-1)^{\tilde a \tilde b}b\circ a$
\item[c)]
{\em  Unit}:

$\forall \in a: 1 \circ a =  a \circ 1 = a$

\item[d)]
{\em Invariance}:

$\eta(a,b\circ c) = \eta(a\circ b,c)$

\end{itemize}
The grading for super Frobenius algebras carries over verbatim.

\subsection{Cobordisms}
\subsubsection{Definition} Let $\mathcal {COB}$ be the category whose objects are
one--dimensional closed oriented (topological) manifolds
considered up to orientation preserving homeomorphism and
whose morphisms are  cobordisms of these objects.

I.e.\ $\Sigma \in Hom(S_{1},S_{2})$ if
$\Sigma$ is an oriented surface with boundary $\del\Sigma \equiv
-S_{1}\coprod S_{2}$.

The composition of morphisms is given by glueing along boundaries with respect to
orientation reversing homeomorphisms.

\subsubsection{Remark}
The operation of disjoint union makes this category into a
monoidal category with unit $\emptyset$.

\subsubsection{Remark} The objects can be chosen to be represented
by disjoint unions of
the circle with the natural orientation $S^1$  and the circle with opposite orientation $\bar S^1$.
Thus a typical object looks like
$S= \coprod_{i \in I}S^{1}\coprod_{j\in J} \bar S^{1}$.
Two standard morphisms are given by the cylinder, and thrice punctured sphere.

\subsubsection{Definition}
Let $\mathcal{VECT}_{k}$ be the monoidal category of finite dimensional
$k$--vector spaces with linear morphisms with the tensor product providing a
monoidal structure with unit $k$.

\subsubsection {Theorem} (Atiyah, Dijkgraaf, Dubrovin) [A,Dij, Du]
There is a 1--1 correspondence between Frobenius
algebras over $k$ and isomorphism classes of covariant functors of monoidal categories from $\mathcal{COB}$
to $\mathcal{VECT}_{k}$, natural with respect to orientation preserving
homeomorphisms of cobordisms and whose value on cylinders $S \times I \in Hom(S,S)$
is the identity.

Under this identification, the Frobenius algebra $A$ is the image of $S^1$, the
multiplication or rather $\mu$
is the image of a thrice punctured sphere and the metric is the image of an annulus.

\section{Orbifold Frobenius algebras}
\subsection{$G$--Frobenius algebras}
\label{orb}

We fix a finite group $G$ and denote its unit element by $e$.

\subsubsection{Definition}
 A {\em G--twisted Frobenius algebra} (FA) over
a field  $k$ of characteristic 0 is
$<G,A,\circ,1,\eta,\varphi,\chi>$, where

\begin{tabular}{ll}
$G$&finite group\\
$A$&finite dim $G$-graded $k$--vector space \\
&$A=\oplus_{g \in G}A_{g}$\\
&$A_{e}$ is called the untwisted sector and \\
&the $A_{g}$ for $g \neq
e$ are called the twisted sectors.\\
$\circ$&a multiplication on $A$ which respects the grading:\\
&$\circ:A_g \otimes A_h \rightarrow A_{gh}$\\
$1$&a fixed element in $A_{e}$--the unit\\
$\eta$&non-degenerate bilinear form\\
&which respects grading i.e. $g|_{A_{g}\otimes A_{h}}=0$ unless
$gh=e$.\\
\end{tabular}

\begin{tabular}{ll}
$\varphi$&an action  of $G$ on $A$ 
(which will be  by algebra automorphisms), \\
&$\varphi\in \mathrm{Hom}(G,\mathrm{Aut}(A))$, s.t.\
$\varphi_{g}(A_{h})\subset A_{ghg^{-1}}$\\
$\chi$&a character $\chi \in \mathrm {Hom}(G,k^{*})$ \\

\end{tabular}

\vskip 0.3cm

\noindent Satisfying the following axioms:

\noindent{\sc Notation:} We use a subscript on an element of $A$ to signify that it has homogeneous group
degree  --e.g.\ $a_g$ means $a_g \in A_g$-- and we write $\varphi_{g}:= \varphi(g)$ and $\chi_{g}:= \chi(g)$.

\begin{itemize}

\item[a)] {\em Associativity}

$(a_{g}\circ a_{h}) \circ a_{k} =a_{g}\circ (a_{h} \circ a_{k})$
\item[b)] {\em Twisted commutativity}

$a_{g}\circ a_{h} = \varphi_{g}(a_{h})\circ a_{g}$
\item[c)]
{\em $G$ Invariant Unit}:

$1 \circ a_{g} = a_{g}\circ 1 = a_g$

and

$\varphi_g(1)=1$
\item[d)]
{\em Invariance of the metric}:

$\eta(a_{g},a_{h}\circ a_{k}) = \eta(a_{g}\circ a_{h},a_{k})$

\item[i)]
{\em Projective self--invariance of the twisted sectors}

$\varphi_{g}|A_{g}=\chi_{g}^{-1}id$

\item[ii)]
{\em $G$--Invariance of the multiplication}

$\varphi_{k}(a_{g}\circ a_{h}) = \varphi_{k}(a_{g})\circ  \varphi_{k}(a_{h})$

\item[iii)]

{\em Projective $G$--invariance of the metric}

$\varphi_{g}^{*}(\eta) = \chi_{g}^{-2}\eta$

\item[iv)]
{\em Projective trace axiom}

$\forall c \in A_{[g,h]}$ and $l_c$ left multiplication by $c$:

$\chi_{h}\mathrm {Tr} (l_c  \varphi_{h}|_{A_{g}})=
\chi_{g^{-1}}\mathrm  {Tr}(  \varphi_{g^{-1}} l_c|_{A_{h}})$
\end{itemize}

An alternate choice of data is given by a one--form $\eps$, the co--unit
with $\eps \in A_e^*$ and a three--tensor
$\mu \in A^* \otimes A^* \otimes A^*$  which is of group degree $e$,
 i.e. $\mu|_{A_{g}\otimes A_{h}\otimes A_{k}}=0$ unless
$ghk=e$.

The relations between $\eta,\circ$ and $\eps,\mu$ are analogous to those of \ref{convert}.

Again, we denote by $\rho\in A_e$ the element dual to $\eps\in A_e^*$ and Poincar\'e dual to $1\in A_e$.

\subsubsection{Remarks}
\begin{itemize}
\item[1)] $A_{e}$ is central by twisted commutativity
and $\la A_{e},\circ,\eta|_{A_{e}\otimes A_e},1\ra$ is a Frobenius algebra.

\item[2)] All $A_{g}$ are $A_{e}$-modules.
\item [3)] Notice that  $\chi$ satisfies the
following equation which completely determines it in terms of $\varphi$.
Setting $h=e,c=1$ in axiom iv)
\begin{equation}
\dim{A_{g}}= \chi_{g^{-1}}{\rm Tr} (\varphi_{g}|_{A_{e}})
\end{equation}
by axiom iii) the action of $\varphi$ on $\rho$ determines $\chi$ up to a sign
\begin{equation}
\chi_{g}^{-2} = \chi_{g}^{-2}  \eta(\rho,1) =
\eta(\varphi_{g}(\rho),\varphi_{g}(1))=  \eta(\varphi_{g}(\rho),1)
\end{equation}
\item[4)] Axiom iv) forces the $\chi$ to be group homomorphisms, so it would be enough to assume in the
data that they are just maps.
\end{itemize}

\subsubsection{Proposition}
The $G$ invariants $A^{G}$ of a $G$--Frobenius algebra $A$
form an associative and commutative algebra with unit.
This algebra with the induced bilinear form
is a Frobenius algebra if and only if
$\sum_{g}\chi_{g}^{-2}=|G|$. If $k={\bf C}$  and $\chi \in \mathrm {Hom}(G,U(1))$ this
implies $\forall g: \chi_g = \pm 1$

\medskip

{\bf Proof.}
Due to axiom ii) the algebra is associative and commutative.
And since $1$ is $G$ invariant, the algebra has a unit.

Now suppose  $\sum_{g}\chi_{g}^{-2} =|G|$.
Then $\eta|_{A^{G}\otimes A^{G}}$ is non--degenerate:
Let $a \in A^{G}$ and choose $b\in A$ s.t. $\eta(a,b)\neq 0$.
Set $\tilde b = \frac{1}{|G|} \sum_{g \in G} \varphi_{g}(b)\in A^G$. Then
\begin{eqnarray*}
\eta(a,\tilde b) &= & \frac{1}{|G|}\sum_{g\in G}
\eta(a,\varphi_{g}(b))= \frac{1}{|G|}\sum_{g\in G}
\eta(\varphi_{g}(a)\varphi_{g}(b))
= \frac{1}{|G|} \sum_{g\in G} \chi_{g}^{-2}\eta(a,b)\\
 &=& \eta(a,b) \neq 0
\end{eqnarray*}

On the other hand if  $\eta|_{A^{G}\otimes A^{G}}$ is non--degenerate
then let
$a,b \in A^{G}$ be such that $\eta(a,b) = 1$. It follows:
$$
1= \eta(a,b) = \frac{1}{|G|} \sum_{g\in G}\eta(\varphi_{g}(a),
\varphi_{g}(b)) = \frac{1}{|G|} \sum_{g\in G} \chi_{g}^{-2}\eta(a,b)
=\frac{1}{|G|} \sum_{g\in G} \chi_{g}^{-2}
$$
so that $\sum_{g}\chi_{g}^{-2} =|G|$.

The last statement follows from the simple fact that since $\forall g
\in G: |\chi_g|=1$ and if $\sum_{g}\chi_{g}^{-2} =|G|$ then
$\chi_g^{-2}=1$ and hence $\chi_g \in \{-1,1\}$.

\subsubsection{Definition}
A $G$--Frobenius algebra is called an {\em orbifold model} if the data
$\la A^G, \circ, 1\ra$ can be augmented by a compatible metric to yield a Frobenius
algebra. In this case, we call the Frobenius algebra $A^G$ a
{\em  $G$--orbifold Frobenius algebra}.

\subsection{Super-grading}

We can enlarge the framework by considering super--algebras rather than
algebras. This will introduce the standard signs.

\subsubsection{Definition} A G-twisted Frobenius super--algebra over
a field  $k$ of characteristic 0 is
$<G,A,\circ,1,\eta,\varphi,\chi>$, where

\begin{tabular}{ll}
$G$&finite group\\
$A$&finite dimensional  $\Z2 \times G$-graded $k$--vector space \\
&$A= A_0 \oplus A_1= \oplus_{g \in G}(A_{g,0}\oplus A_{g,1}) = \oplus_{g \in G}A_{g}$\\
&$A_{e}$ is called the untwisted sector and is even. \\
&The $A_{g}$ for $g \neq  e$ are called the twisted sectors.\\
$\circ$&a multiplication on $A$ which respects both gradings:\\
&$\circ:A_{g,i} \otimes A_{h,j} \rightarrow A_{gh,i+j}$\\
$1$&a fixed element in $A_{e}$--the unit\\
$\eta$&non-degenerate even bilinear form\\
&which respects grading i.e. $g|_{A_{g}\otimes A_{h}}=0$ unless
$gh=e$.\\
$\varphi$&an action by even algebra automorphisms  of $G$ on $A$, \\
&$\varphi\in \mathrm{Hom}_{k-alg}(G,A)$, s.t.\
$\varphi_{g}(A_{h})\subset A_{ghg^{-1}}$\\
$\chi$&a character $\chi \in \mathrm {Hom}(G,k^{*})$ \\
&or if $k = {\bf C}$, $\chi \in \mathrm {Hom}(G,U(1))$
\end{tabular}

\vskip 0.3cm

\noindent Satisfying the  axioms a)--d) and i)--iii) of a $G$--Frobenius
algebra with the
following alteration:

\begin{itemize}

\item[b$^{\sigma}$)] {\em Twisted super--commutativity}

$a_{g}\circ a_{h} = (-1)^{\tilde a_g\tilde a_h} \varphi_{g}(a_{h})\circ a_{g}$

\item[iv$^{\sigma}$)]
{\em Projective super--trace axiom}

$\forall c \in A_{[g,h]}$ and $l_c$ left multiplication by $c$:

$\chi_{h}\mathrm {STr} (l_c  \varphi_{h}|_{A_{g}})=
\chi_{g^{-1}}\mathrm  {STr}(  \varphi_{g^{-1}} l_c|_{A_{h}})$
\end{itemize}
where $\mathrm{STr}$ is the super--trace.

\subsubsection{Operations}

{\sc Restriction.}

If $H \subset G$ and $A= \oplus_{g\in G}A_g$ then $\tilde A:= \bigoplus_{h\in H} A_h$ is naturally a
 $H$--Frobenius algebra.

{\sc Direct Sum.}
Given a $G$--Frobenius algebra $A$ and an $H$--Frobenius algebra $B$ then $A\oplus B$
is  naturally a $G\times H$--Frobenius algebra with the graded pieces $(A\oplus B)_{(g,h)}= A_g \oplus B_h$.

We define the direct sum of two $G$--Frobenius algebras to be the $G$--Frobenius subalgebra corresponding to the diagonal
$\Delta:G \rightarrow G\times G$ in $A\oplus A$.

{\sc Tensor product.}
Given a $G$--Frobenius algebra $A$ and an $H$--Frobenius algebra $B$ then $\bigoplus_{(g,h)} (A_g\otimes B_h)$
is  naturally a $G\times H$--Frobenius algebra $(A\otimes B)_{(g,h)}= A_g \otimes B_h$.

We define the tensor product of two $G$--Frobenius algebras to be the
$G$--Frobenius subalgebra corresponding to the diagonal $G \rightarrow G\times G$ in
$A \otimes A$.

{\sc Braided Tensor product}.
If $A$ and $B$  are two $G$--Frobenius algebras with the same character $\chi$, we can
define a braided tensor product structure
on $A\otimes B$ by setting $(A\otimes B)_g := \oplus_{k\in G} A_k \otimes B_{k^{-1}g}$.
For the multiplication we use the sequence

\begin{multline}
A_k \otimes B_{k^{-1}g} \otimes A_l \otimes B_{l^{-1}h}
\stackrel{(id \otimes id \otimes  \varphi_{k^{-1}g}\otimes id)\circ \tau_{2,3}}{\longrightarrow}
A_k \otimes A_{k^{-1}glg^{-1}k}\otimes B_{k^{-1}g }\otimes  B_{l^{-1}h }\\
\stackrel {\circ\otimes \circ}{\rightarrow} A_{glg^{-1}k} \otimes B_{k^{-1}gl^{-1}h}
\end{multline}

and $\bigoplus_k (\varphi_k \otimes \varphi_ {kh^{-1}})$ for the action of $h$ on
$\bigoplus_k A_k \otimes B_{k^{-1}g} $

\subsubsection{Remark}
If one thinks in terms of cohomology of spaces the direct sum corresponds to the disjoint union
and the tensor product corresponds to the Cartesian product. The origin of the braided tensor product,
however, is not clear yet.

\subsection{Geometric model -- spectral flow}
The axioms of the $G$--Frobenius algebra are well suited for taking the quotient, since the
 multiplication is
$G$--invariant. However, this is not the right framework for a
geometric interpretation. In order to accommodate a more natural
coboundary description, we need the following definition which
corresponds to the physical notion of Ramond ground states:

\subsubsection{Definition}
A {\em Ramond $G$--algebra} over
a field  $k$ of characteristic 0 is
$<G,V,\bar\circ,v,\bar\eta,\bar\varphi,\chi>$

\begin{tabular}{ll}
$G$&finite group\\
$V$&finite dim $G$-graded $k$--vector space \\
&$V=\oplus_{g \in G}V_{g}$\\
&$V_{e}$ is called the untwisted sector and \\
&the $V_{g}$ for $g \neq
e$ are called the twisted sectors.\\
$\bar\circ$&a multiplication on $V$ which respects the grading:\\
&$\bar\circ:V_g \otimes V_h \rightarrow V_{gh}$\\
$v$&a fixed element in $V_{e}$--the unit\\
\end{tabular}

\begin{tabular}{ll}
$\bar\eta$&non-degenerate bilinear form\\
&which respects grading i.e. $\bar\eta|_{V_{g}\otimes V_{h}}=0$ unless
$gh=e$.\\
$\bar\varphi$&an action by  of $G$ on $V$, \\
&$\bar\varphi\in \mathrm{Hom}(G,\mathrm{Aut}(A))$, s.t.\
$\bar\varphi_{g}(V_{h})\subset V_{ghg^{-1}}$\\
$\chi$&a character $\chi \in \mathrm {Hom}(G,k^{*})$ \\
\end{tabular}

\vskip 0.3cm

\noindent Satisfying the following axioms:

\noindent{\sc Notation:} We use a subscript on an element of $V$ to
 signify that it has homogeneous group
degree  --e.g.\ $v_g$ means $v_g \in V_g$-- and we write $\bar\varphi_{g}:= \bar\varphi(g)$ and $\chi_{g}:= \chi(g)$.

\begin{itemize}

\item[a)] {\em Associativity}

$(v_{g}\bar\circ v_{h}) \bar\circ v_{k} =v_{g}\bar\circ (v_{h} \bar\circ v_{k})$
\item[b')] {\em Projective twisted commutativity}

$v_{g}\bar\circ v_{h} = \chi_{g} \bar\varphi_{g}(v_{h}) \bar\circ v_{g}
= \bar\varphi_{g}(v_{h} \bar\circ v_{g}) $
\item[c')]
{\em Projectively invariant unit}:

$v \bar\circ v_{g} = v_{g}\bar\circ v = v_g$

and

$\bar\varphi_g(v)=\chi_g v$

\item[d)]
{\em Invariance of the metric}:

$\eta(v_{g},v_{h}\bar\circ v_{k}) = \eta(v_{g}\bar\circ v_{h},v_{k})$
\item[1')]
{\em Self--invariance of the twisted sectors}

$\bar\varphi_{g}|V_{g}=id$

\item[2')]
{\em Projective $G$--invariance of multiplication}

$\bar\varphi_{k}(v_{g}\bar\circ v_{h}) = \chi_{k}\bar\varphi_{k}(v_{g})
\bar\circ  \bar\varphi_{k}(v_{h})$

\item[3')]

{\em  $G$--Invariance of metric}

$\bar\varphi_{g}^{*}(\bar \eta) = \bar \eta$

\item[4')]
{\em Trace axiom}

$\forall c \in V_{[g,h]}$ and $l_c$ left multiplication by $c$: \\
$\mathrm {Tr} (l_c \circ \bar\varphi_{g}|_{V_{h}})= \mathrm {Tr}
(\bar\varphi_{h^{-1}}\circ l_c|_{V_{g}})$
\end{itemize}

\subsubsection{Definition}
A {\em state--space} for a $G$--Frobenius algebra $A$
 is a quadruple $\la V,v, \bar\eta,\bar\varphi\ra$, where

\begin{tabular}{ll}
$V$& is a $G$-graded free rank one $A$--module: $V= \bigoplus_{g\in G}V_g$,\\
$v$&is a fixed generator of $V$-- called the vacuum,\\
$\bar\eta$& is non--degenerate bilinear form on $V$\\
$\bar\varphi$& is a linear $G$-action on $V$  fixing $v$ projectively\\
&i.e. $\bar\varphi(g)(span(v)) \subset span(v)$ \\
\end{tabular}

such that these structures are compatible with those of $A$ (we denote $\bar\varphi_g:=\bar\varphi(g)$)
\begin{itemize}
\item[a)] The action of $A$ respects the grading: $A_g V_h \subset V_{gh}$.
\item[b)] $V_{h}$ is a rank one free $A_h$--module and
$V_{h}= A_{h}v$.
\item[c)] $\bar\varphi_{g}(av) = \varphi_{g}(a) \bar\varphi(v): \; \forall a \in
A, v \in V$
\item[d)] For $a,b \in A:
\bar\eta(av,bv) = \eta(a,b)$
\item [e)] $\forall g, h \in G, c \in A_{[g,h]}$
$\forall c \in V_{[g,h]}$ and $l_c$ left multiplication by $c$: \\
$\mathrm {Tr} (l_c \circ \bar\varphi_{g}|_{V_{h}})= \mathrm {Tr}
 (\bar\varphi_{h^{-1}}\circ l_c|_{V_{g}})$
\end{itemize}

\subsubsection{Definition}
We call two state spaces isomorphic, if there is an $A$--module isomorphism between the two.

Since state spaces are free rank one $A$--modules,
it is clear that all automorphisms are rescalings
of $v$.

\subsection{Proposition} Given a $G$--twisted Frobenius algebra $A$
there is  a unique state space up to
isomorphism and the form $\bar\eta$ is
$G$--invariant (i.e.\ $\bar\varphi_g^*(\bar\eta) = \bar\eta$).

{\bf Proof.} We start with a free rank one  $A$--module $V$ and reconstruct all other data.
The $G$--grading on $V$ is uniquely determined from that of $A$ by axiom b) and this grading
satisfies axiom a).
Up to isomorphism, we may assume a generator $v\in V$ is fixed, then
axiom d) determines $\bar\eta$ from  $\eta$.
We denote the Eigenvalue of $\bar\varphi_g$ on $v$ by $\l_g$: $\bar\varphi_g v = \l_g v$.
Notice that due to c) $\bar\varphi$ is determined by $\l_g$. Using axioms b), c) and e),
we find that $ \l_{g}=\chi_g$, thus fixing the $G$--action $\bar\varphi$.

Namely with $c=1$ and $h=e$ in e):
$$
{\rm Tr} ( \bar\varphi_{g}|_{V_{e}}) =
\l_{g}{\rm Tr} (\varphi_{g}|_{A_{e}})
= \l_{g} \chi_g^{-1} {\rm Tr}
(\varphi_e|_{A_{g}})
=\l_{g} \chi_g^{-1} {\rm Tr}
(\varphi_e|_{V_{g}})
$$
The equality $\bar\varphi_{g}=\varphi_g\l_{g}$ implies that $\bar\eta$ is
$\bar\varphi$ invariant:
$$
\bar\eta(\bar\varphi_{g}(av),\bar\varphi_{g}(bv)) =
\l_{g}^{2}\eta(\varphi_{g}(a),\varphi_{g}(b)) =
\l_{g}^{2} \chi_{g}^{-2}\eta(a,b)= \bar\eta(av,bv)
$$
In general:
\begin{eqnarray*}
\mathrm {Tr} (l_c \circ \bar\varphi_{g}|_{V_{h}})&=&
\l_{g}\mathrm {Tr} (l_c \circ  \varphi_{g}|_{A_{h}})\\
&=&\l_{g} \chi_g^{-1} \chi_{h^{-1}} \mathrm {Tr}  (\varphi_{h^{-1}}\circ l_c|_{A_{g}})\\
&=&\l_{g} \chi_g^{-1} \chi_{h^{-1}} \l_{h^{-1}}^{-1}
\mathrm {Tr}  (\bar\varphi_{h^{-1}}\circ l_c|_{V_{g}})\\
&=& \mathrm {Tr}  (\bar\varphi_{h^{-1}}\circ l_c|_{V_{g}})
\end{eqnarray*}
so that with this choice of $\bar\varphi$ axiom e) is satisfied.

\subsubsection{Remark} A state space inherits an associative multiplication $\bar\circ$
 with unit from $A$ via
\begin{equation}
\forall a,b \in A:  (av)\bar\circ (bv) := (a\circ b)v
\end{equation}
This multiplication makes it into a $G$--Ramond algebra.

This fact leads us to the following definitions:
\subsubsection{Definition}
The {\em Ramond space} of a $G$--Frobenius algebra $A$ is the
state--space given by the $G$--graded vector--space
$$V:= \oplus_{g}V_{g}:=\bigoplus_{g}A_{g}\otimes k$$
together with the $G$-action $\bar\varphi := \varphi \otimes
\chi$, the induced metric $\bar g$ and the induced multiplication
$\bar\circ$ and fixed element $v=1\otimes 1$.

\subsubsection{Theorem}
\label{Ramond}
There is a one--to--one correspondence between  isomorphism classes
of $G$--Ramond algebras and
$G$-Frobenius algebras.

{\bf Proof.} The association of a Ramond space to a $G$--Frobenius algebra provides the correspondence.
The inverse being the obvious reverse twist by $\chi$.

\subsubsection{Remark}
In the theory of singularities, the untwisted  sector of the
Ramond space corresponds to the forms $H^{n-1}(V_{\eps},{\bf C})$
while the untwisted sector of the $G$--twisted Frobenius algebra
corresponds to the Milnor ring [Wa]. These are naturally
isomorphic, but have different $G$--module structures. In that
situation, one takes the invariants of the Ramond sector, while we
will be interested in invariants of the $G$--Frobenius algebra and
not only the untwisted sector (cf.\ [K5] and see also \S 7).

\section{Bundle cobordisms, finite gauge groups, orbifolding and $G$--Ramond algebras}

In this section, we introduce two cobordism categories which
correspond to $G$--orbifold Frobenius algebras and Ramond
$G$--algebras, respectively.  Again $G$ is a fixed finite group.

\subsection{Bundle cobordisms}

In all situations, gluing along boundaries will 
induce the composition and the disjoint
union will provide a monoidal structure.

\subsubsection{Definition} Let $\mathcal {GBCOB}$ be the category
whose objects are principal $G$--bundles over one--dimensional
closed oriented (topological) manifolds, pointed over each
component of the base space, whose morphisms are cobordisms of
these objects (i.e.\ principal $G$--bundles over oriented surfaces
with pointed boundary).

More precisely, $B_{\Sigma} \in Hom(B_{1},B_{2})$ if
$\Sigma$ is an oriented surface
with boundary $\del\Sigma =
-S_{1}\coprod S_{2}$
and $B_{\Sigma}$ is
a bundle on $\Sigma$ which restricts to $B_{1}$ and $B_{2}$ on the
boundary.

The composition of morphisms is given by gluing along
boundaries with respect to
orientation reversing homeomorphisms
on the base and covering bundle isomorphisms which align the base--points.

\subsubsection{Remark}
The operation of disjoint union makes this category into monoidal category with unit $\emptyset$
formally regarded
as a principal $G$  bundle over  $\emptyset$.

\subsubsection{Remark}
Typical objects are bundles $B$ over \\
$S= \coprod_{i \in I}S^{1}\coprod_{j\in J} \bar S^{1}$.

Let $(S^{1},\nu)$ $\nu \in S^1)$ be a pointed $S^1$.

\subsubsection{Structure Lemma}
\label{structure}
The space $Bun(S^{1},G)$ of $G$ bundles on $(S^{1},\nu)$ can be described as follows:
$$Bun(S^{1},G) = (G \times F)/G$$
where  $F$ is a generic fibre regarded as a principal $G$-space
and $G$ acts on itself by conjugation and the quotient is taken by
the diagonal action.  The space $Bun^*(S^{1},G)$ of pointed $G$
bundles on $S^{1}$ is a $G$ bundle over $Bun(S^{1},G)$.
Evaluating the monodromy give a projection of $Bun^*(S^{1},G)$ onto $G$ whose 
fiber over $g \in G$ are the centralizers $Z(g)$.

\subsubsection{Remark} Usually one uses the identification
$Bun(M,G)=\mathrm{Hom}(\pi_1(M),G)/G $, which we also use in the
proof. However, for certain aspects of the theory -- more
precisely to glues and to include non--trivial characters-- it is
vital to include a point in the bundle and a trivialization rather
then just a point in the base.

{\bf Proof of the Structure Lemma.} Given a pointed principal $G$
bundle $(B,S^1,\pi$$,F,G)$,
 $b\in B$ we set $\nu = \pi(b) \in S^1$. The choice of $b \in F =B_{\nu}$ gives an identification of $F$ with $G$, i.e.\ we
let $\b:G\mapsto F$ be the admissible map in the sense of [St]
that satisfies $\b(e)=b$. We set $g\in G$ to be the element
corresponding to the monodromy around the generator of
$\pi_1(S^1)$. Notice that since we fixed an admissible map
everything is rigid -- there are no automorphisms-- and the
monodromy is given by an element, not a conjugacy class. Thus we
associate to $(B,b)$ the tuple $(g,b)$.

Vice versa, given $(g,b)$ we start with the pointed space
$(S^1,\nu)$ and construct the
 bundle with fiber $G$,
monodromy $g$ and marked point $b=\b(e)\in B_{\nu}$.

The bijectiveness of this construction follows by the classical results quoted above [St].

The choice in this construction corresponds to a choice of a point
$b \in B$. Changing $b$ amounts to changing $\b$ and the monodromy
$g$. Moreover,  moving the point $\nu= \pi(b)$ and moving $b$
inside the fiber by parallel transport keeps everything fixed.
Moving $b$ inside the fixed fiber by translation  (once $\nu$ is
fixed)  corresponds to translation by the group action in the
fiber i.e.\ the translation action of $G$ on $F$ and simultaneous
conjugation of the monodromy. Hence, the first claim follows.
The second claim follows from the fact that we are dealing with pricipal 
bundles. Finally evaluating the monodromy of two pointed 
bundles yields the same result if the bundles are isomorphic and the
points are related by a shift in the centralizer of the monodromy.

This observation leads us to the following definition:

\subsubsection{Definition}
We call a bundle over a closed one--dimensional manifold
rigidified if its components are labelled and  the bundle is
pointed above each component of the base and  trivializations
around the projection of the marked points to the base are fixed.
We denote such a bundle $(B,b)$ where $b \subset B: b=
(b_0,\dots,b_n)$ is the set of base--points for each component.

Furthermore, if a given surface has genus zero we realize it in
the plane as a pointed disc with all boundaries being $S^1$. We
label the outside circle by $0$.

 If  $\pi$ is the bundle projection,  we set $x_i:= \pi(b_i)$,
 we call $b_0$ the base--point,
the $B_{x_i}$ the special fibers and $B_{x_0}$ the initial fiber.

\subsection{Rigidified bundle cobordisms and $G$--Frobenius algebras.}

\subsubsection{Definition} Let $\mathcal {GBCOB}^*$ be the category whose objects are
rigidified
principal $G$--bundles over one--dimensional closed oriented (topological) manifolds
considered up to pull--back under orientation preserving homeomorphism of the base respecting the markings
 and
whose morphisms are  cobordisms of these objects (i.e. principal
$G$--bundles over oriented surfaces with boundary together with
rigidification on the boundary, a choice of base--point $x_0\in
\del\Sigma$ compatible with the rigidification and a choice of
curves $\Gamma_i$ from $x_0$ to $x_i$ which identify the
trivializations via parallel transport, where we used the notation
above.).

I.e.\ objects are bundles $B$ over $S= \coprod_{i \in
I}S^{1}\coprod_{j\in J} \bar S^{1}$ with base--points  on the
various components $\bar b_1=(b_{i}\in B|_{S^1_i}:i\in I), \bar
b_2=(b_j\in B|_{\bar S^1_j}:j\in J)$ and $B_{\Sigma} \in
Hom((B_{1},\bar b_1),(B_{2},\bar b_2))$ if $\Sigma$ is an oriented
surface with boundary considered up to orientation preserving
homeomorphism with boundary $\del\Sigma = -S_{1}\coprod S_{2}$ --
again up to homeomorphism -- and $B_{\Sigma}$ is a bundle on
$\Sigma$ which restricts to $B_{1}$ and $B_{2}$ on the boundary
together with rigidification data for the boundary i.e. $(\bar
b_1,\bar b_2)$ with $\bar b_1 \subset B_1$ and $\bar b_2 \subset
B_2$. We will call $S_1$ the inputs and $S_2$ the outputs.

The extra structure of curves and base--point allows us to identify the special fibres
with the initial fibre
via parallel transport. Thus we can describe the rigidification data in terms of $b_0$,
and the
group elements $g_i$ defined via $\Gamma_{i}(b_0) = g_{i}b_i$ .

The composition of morphisms is given by gluing along boundaries with respect to
orientation reversing homeomorphisms of the base and covering  bundle isomorphisms
identifying the base--points.
\subsubsection{Remark}
The operation of disjoint union makes this category into a
monoidal category with unit $G$ regarded
as a bundle over  $\emptyset$ with base point $e$.

\subsubsection{Construction}
\label{Construction}
We define the pointed bundle $(g,h)$ to be the bundle wich is obtained
by glueing $I \times G$ via the identification $(0,e)\sim(1,g)$ and
marking the point $(0,h)$,
where $I=[0,1]$ the standard interval. 

This produces all monoidal generators.
 
\subsubsection{Remark}
By the Construction \ref{Construction}  up to reversing the orientation,
we can produce the monoidal generators of
objects of $\mathcal {GBCOB}*$ with the objects coming from $G\times G$.

A generating object can thus be depicted by an oriented circle
with labels $(g,h)$ where $(g,h) \in G\times G$. We use the
notation $(g,h)$ for positively oriented circles and
$\overline{(g,h)}$ for negatively oriented circles. We will
consider functors $V$ with an involutive property we have that we
will have $V(\overline{(g,h)})\simeq V((g^{-1},h))^* $ where $*$
denotes the dual. General objects are then just disjoint unions of
these, i.e.\ tuples $(g_i,h_i)$. Homomorphisms on the generators
are given by the trivial bundle cylinder with different
trivializations on both ends represents the natural diagonal
action of $G$ by conjugation and translation described in
$\ref{structure}$, so that this diagonal $G$--action is realized
in terms of  cobordisms.

The natural $G\times G$ action, however, cannot be realized by these cobordisms and we
would like to enrich our situation to this case by adding morphisms corresponding to the
$G$ action on the trivializations.

\subsubsection{Definition} Let $\mathcal {RGBCOB}$
be the category obtained from $\mathcal {GBCOB*}$ by adding the
following morphisms. For any n--tuple $(k_1,\dots,k_n): k_i \in G$
and any object $ (g_i,h_i):i=1,\dots n$ we set

$$\tau(k_1,\dots,k_n)(g_i,h_i)_{i\in \{1,\dots n\}}:=(g_i,k_ih_i)_{i\in \{1,\dots n\}}$$

We call these morphisms of type II and the morphisms coming from
$\mathcal GBCOB*$ of type I. We also sometimes write $II_k$ for
$\tau_k$.

\subsubsection{Remark}
There is a natural forgetful functor from $\mathcal {RGBCOB}$
to $\mathcal {GBCOB}$.

Given a character $\chi\in Hom (G,k^*)$ we can form the fibre--product $G\times_{\chi}k$.
This gives a functor $k[G]_\chi$ from the category $\mathcal {VECT}_k$ to $k[G]-\mathcal{MOD}$,
the category of $k[G]$--modules.

\subsubsection{Definition}
A $G_{\chi}$--orbifold theory is a monoidal
functor $V$ from $\mathcal{RGBCOB}$ to $k[G]-\mathcal{MOD}$
satisfying the following conditions:
\begin{itemize}
\item[i)] The image of $V$ lies in the image of $k[G]_\chi$.
\item[ii)] Objects of $\mathcal {RGBCOB}$
which differ by  morphisms of type II are mapped to the same object in $k[G]-\mathcal{MOD}$.
\item[iii)]  The value on morphisms of type I does not depend on the choice of
connecting curves and associated choice of trivializations or
base--point.
\item[iv)] The morphisms of type II are mapped to the $G$--action by $\chi$.
\item[v)] The functor is natural with respect to morphisms of the type $\tau(k,\dots,k)$.
I.e.\  $V(\tau(k,\dots,k) \circ \Sigma) =
V(\tau_{out}(k^{-1},\dots,k^{-1}))\circ V(\Sigma) \circ
V(\tau_{in}(k,\dots, k))$; where $\tau_{in}$ and $\tau_{out}$
operator on the inputs and outputs respectively.
\item[vi)]  $V$ associates $id$ to  cylinders
$B \times I,(b,0)\in B\times 0, (b,1)\in B\times 1$ considered
as cobordisms from $(B,b)$ to itself.
\item[vii)] $V$ is involutive: $V(\bar S)= V^*(S)$ where $*$ denotes the dual vector space
with induced k[G]--module structure. In accordance, the morphism
of type II  commute with involution, i.e.\ they are mapped to the
$G$--action by $\chi^{-1}$.
\item[viii)] The functor is natural with respect to orientation preserving homeomorphisms
of the underlying surface of a cobordism
and pull--back of the bundle.
\end{itemize}

\subsubsection{Corollary}
\label{homotopy}
$G_{\chi}$--orbifold theories are homotopy invariant.

{\bf Proof.} By standard arguments using naturality and vi) a homotopy of an object $S$ induces an identity
on its image. More precisely given a homotopy of objects $f_t: S\mapsto S$ , it induces a map
$F:S\times I \mapsto S\times I$ and we have the commutative diagram:

\begin{tabular}{ccc}
$V(S)$&$\stackrel{V(S\times I)=id}{\rightarrow}$&$V(S)$\\
$V(f_0) \downarrow$&&$\downarrow V(f_1)$\\
$V(S)$&$\stackrel{V(S\times I)=id}{\rightarrow}$&$V(S)$\\
\end{tabular}

In particular, if $V(f_0)= id$ then  $V(f_1)=id$.

\subsubsection{Remark}

Given a choice of connecting curves,
we can identify all fibres over special points.
Therefore after fixing one identification of a fibre with $G$, we can identify the other marked points as
translations of points of parallel transport and label them by group elements, which we will do.

The action of $\tau(k,\dots , k)$ corresponds to a change of
identification for one point and simultaneous re--gauging of all
other points via this translation, i.e.\ a diagonal gauging.
Therefore given a cobordism, we can fix an identification of all
fibres with $G$.

\subsubsection{Definition and Notation}

We will fix some standard bundle cobordisms pictured below.

\medskip

\epsfxsize = \textwidth
\epsfbox{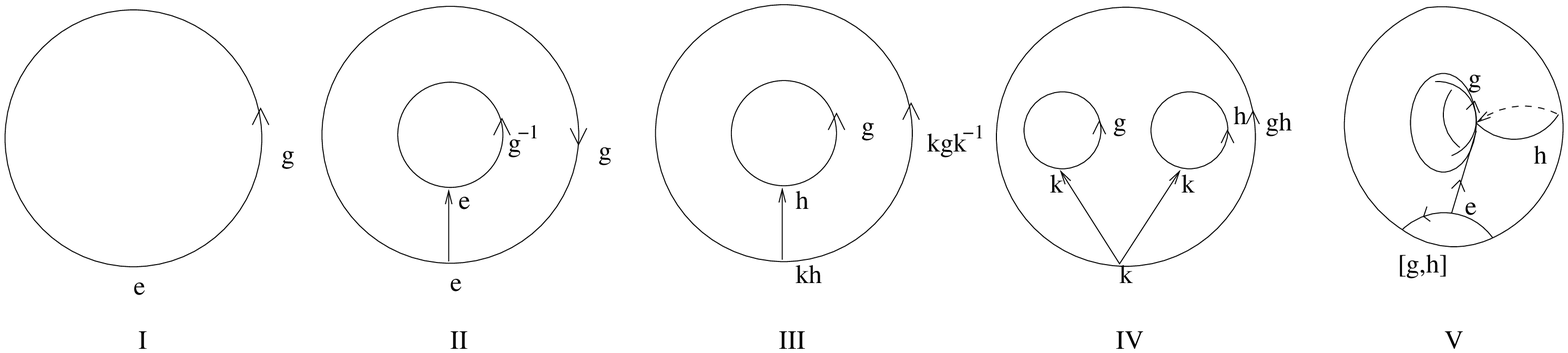}

\medskip

I: The {\em standard disc  bundle} is the disc with a trivial bundle and positively oriented boundary
considered as a cobordism between $\emptyset$ and $(e,e)$; it will be denoted by $D$.

II: The {\em standard $g$--cylinder bundle} is the cylinder $S^1\times I$ with the bundle having
monodromy $g$
around $((S^1,0))$ considered as a cobordism between $(g,e) \coprod \overline {(g,e)}$  and $\emptyset$;
it will be denoted by $C_g$.

III: The {\em $(g,h)^k$--cylinder bundle} is the cylinder with a
bundle having monodromy $g$ around $(S^1,0)$ considered as a
cobordism between $(g,h)$ and $(kgk^{-1},kh)$; it will be denoted
by $C^h_{g,k}$.

IV: The {\em standard $(g,h)^k$--trinion bundle}
 is the trinion with the bundle having monodromies $g$ around the first $S^1$ and
$h$ around the second $S^1$ and translations $e$ for $\tau_{01},
\tau_{02}$ considered as a cobordism between $(g,k),(h,k)$ and
$(gh,k)$; it will be denoted by $T^k_{g,h}$.

V: The {\em $(g,h)$--torus bundle} is the once--punctured torus with  the principal $G$--bundle having
monodromies $g$ and $h$ around the two standard cycles
considered as a cobordism between $([g,h],e)$ and $\emptyset$; it will be denoted by $E^k_{g,h}$.

\subsubsection{Lemma}
The bundles over a cylinder that form cobordisms between  $(g,e)$ and $(h,k)$ are
parameterized by  $G$; it is necessary that $g=khk^{-1}$.
These cobordisms are given by the $C_{g,h}^e$.

{\bf Proof.} Given such a bundle over the cylinder $\Sigma_0:=S^1 \times I$ the translation
from $B_{\nu,0}$ to $B_{\nu,1}$ along $\gamma(t) :=(\nu,t)$ is a complete
invariant. We fix this
element $k\in G$.
Since in $\pi_1 (\Sigma_0,(\nu,0))$ $C_1 = \gamma C_2 \gamma^{-1}$
we must have $g=khk^{-1}$. 

\subsubsection{Proposition}
\label{fix}
To fix a $G_{\chi}$--orbifold theory on the objects of the type $(g,e)$ and
to fix a $G_{\chi}$--orbifold theory on the morphisms it suffices to
fix it on
bundles over the standard cylinder bundle $C$, the $(g,h)$--cylinder bundles $C_{g,h}$,
and on the standard trinion $T$.

{\bf Proof.} The first claim follows from the condition ii) and the monoidal structure.
For the second claim first notice that due to the homotopy Lemma \ref{homotopy} $V$ is fixed on $D$.
Also,  any bundle over the cylinder is trivial, furthermore by
v) we may regard the cylinder as a cobordism from $B|_{S^1}$ to itself and after applying
a morphism of type II we can assume that the boundary objects are of the type $(g,e),(h,k)$.
Therefore we know the functor $V$ on all bundles over cylinders.
Dualizing and gluing on cylinders, we find that once the functor is defined on the standard trinion,
it is defined on all bundles over all trinions. Lastly, given any surface, we can choose a decomposition
by a marking into discs, cylinders and trinions, then  gluing determines the value of $V$ on this surface.

\subsubsection{Proposition and Definition}
\label{propdef}
For any  $G_{\chi}$--orbifold theory $V$ set

\begin{tabular}{lcl}
$V((g,e))$&=&$V_g$\\
$V(T^e_{g,h})$&=&$\bar\circ: V_g \otimes V_h \rightarrow V_{gh}$\\
$V(D)(1_k)$&=&$v$\\
$V(C_g)$&=&$\bar\eta|_{V_g\otimes V_{g^{-1}}}$\\
$V(C_{g,h})$&=&$\bar\varphi_h:V_g\mapsto V_{hgh^{-1}}$
\end{tabular}

then this data together with $g$ and $\chi$ form
a Ramond--$G$ algebra $\la G,V,\bar\circ,v,\bar\eta,\bar\varphi,\chi\ra$
which we call the {\em associated $G$--Ramond algebra to $V$}.

{\bf Proof.}
It is clear by \ref{fix} that given a $G_{\chi}$ orbifold theory it is completely fixed by its
associated $G$--Ramond algebra.
The converse is also true:

The axiom a) follows from the standard gluing procedures of TFT. I.e.\
the usual gluing of a disc with 3
 holes from two discs with two holes in two different ways.

\medskip

\epsfxsize =\textwidth
\epsfbox{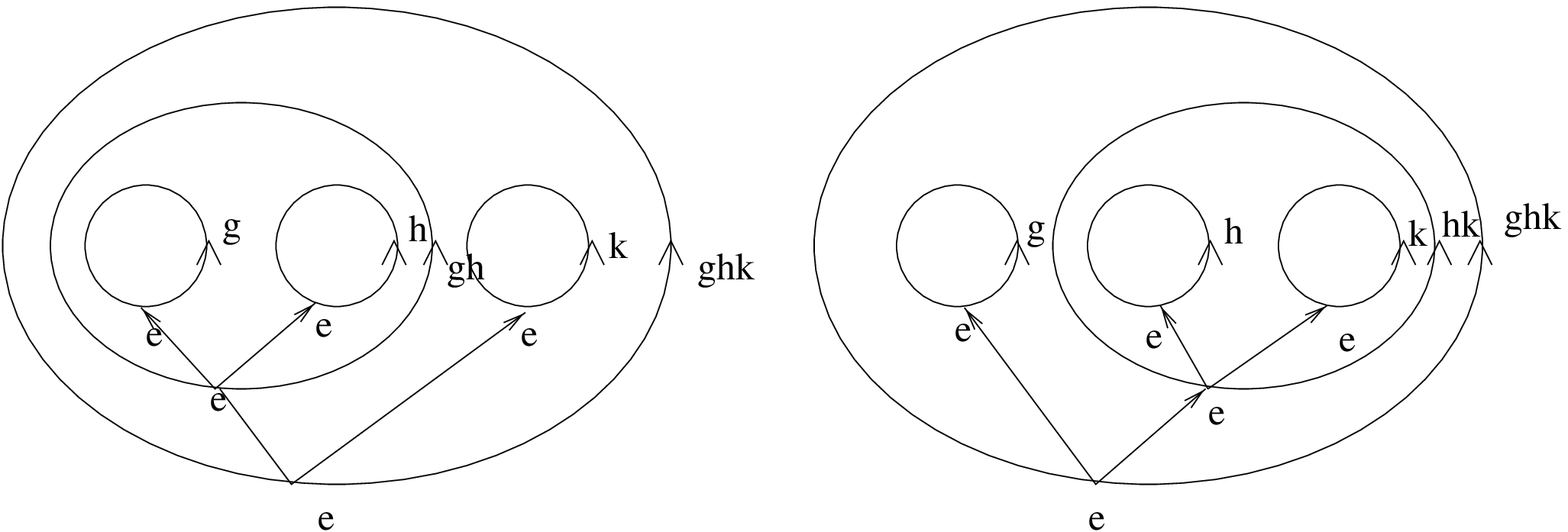}

\medskip

For axiom b') we regard the following commutative diagrams

\begin{tabular}
{@{\extracolsep{0pt}}c@{\extracolsep{0pt}}c@{\extracolsep{0pt}}c@{\extracolsep{0pt}}c@{\extracolsep{0pt}}c@{\extracolsep{0pt}}c@{\extracolsep{0pt}}c@{\extracolsep{0pt}}}
$(g,e)\coprod (h,e)$&$\stackrel{\tau_{12}}{\longrightarrow}$&$(h,e)\coprod (g,e)$&
&$(h,e)\coprod (g,e)$&
&$(g,e)\coprod (h,e)$\\
$\downarrow T_{g,h}^e$&$\Rightarrow$&$\downarrow II$&$\Rightarrow$&$\downarrow III$&
$\Rightarrow$&$\downarrow IV$\\
$(gh,e)$&$\stackrel{id}{\longrightarrow}$&$(gh,e)$
&&$(gh,g)$
&&$(gh,g)$\\
&&&$\Downarrow V$&&&\\
$V_g\otimes V_h$&$\stackrel{\tau_{12}}{\longrightarrow}$&$V_h\otimes V_g$
&$\stackrel{id}{\longrightarrow}$&$V_h\otimes V_g$
&$\stackrel{id}{\longrightarrow}$&$V_h\otimes V_g$\\
$\downarrow \bar\circ$&&&&&&$\downarrow
\bar\varphi_g \bar\circ$\\
$V_{gh}$&$\stackrel{id}{\longrightarrow}$&$V_{gh}$
&$\stackrel{id}{\longrightarrow}$&$V_{gh}$
&$\stackrel{id}{\longrightarrow}$&$V_{gh}$\\
\end{tabular}

where we have used axiom viii) for the first move and axiom iii)
for the second and gluing for the last one.
\medskip

\epsfxsize =\textwidth
\epsfbox{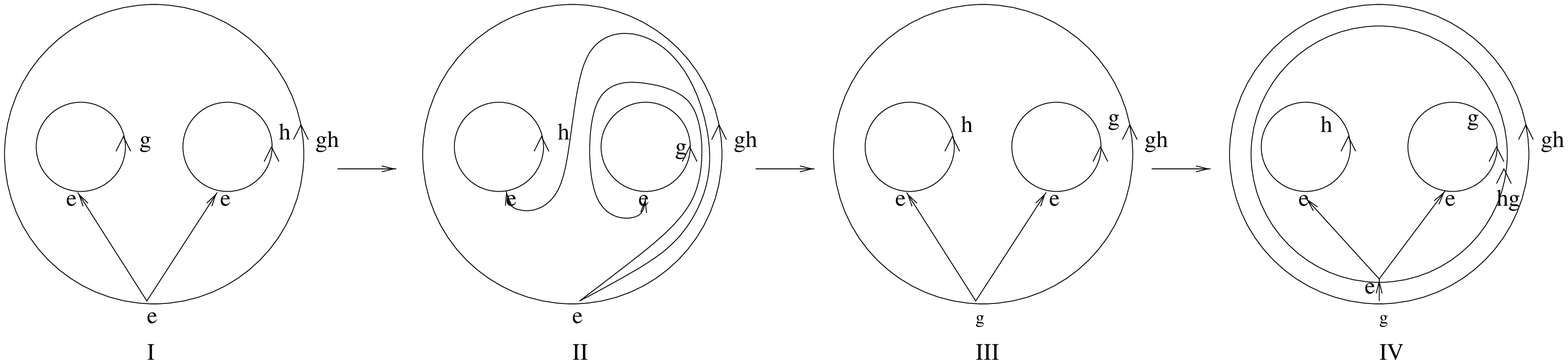}

\medskip

The unit of axiom c) is given by $\bar D^e$. The projective invariance follows from

\begin{tabular}
{@{\extracolsep{0pt}}c@{\extracolsep{0pt}}c@{\extracolsep{0pt}}c@{\extracolsep{0pt}}c@{\extracolsep{0pt}}c@{\extracolsep{0pt}}c@{\extracolsep{0pt}}c@{\extracolsep{0pt}}}
$\emptyset$&\multicolumn{5}{c}{$\stackrel{D}{\longrightarrow}$}&$(e,e)$\\
$\downarrow \|$&&&&&&$\downarrow \| $\\
$\emptyset$&$\stackrel{D}{\rightarrow}$&$(e,e)$&
$\stackrel{C_{(e,e)}^k}{\rightarrow}$&$(e,k)$&
$\stackrel{II_{k^{-1}}}{\longrightarrow}$&$(e,e)$\\
&&&$\Downarrow V$&&&\\
$k$&$\stackrel{v}{\rightarrow}$&$V_e$&$\stackrel{\bar\varphi_k}{\rightarrow}$&$V_e$
&$\stackrel{\bar\chi_{k^{-1}}}{\rightarrow}$&$V_e$\\
$\downarrow \|$&&&&&&$\downarrow \| $\\
$k$&\multicolumn{5}{c}{$\stackrel{1^e}{\longrightarrow}$}&$V_e$\\
\end{tabular}

where in the third line $1_k \mapsto v \mapsto \chi_{k}v\mapsto v$.

The axiom  d) of the invariance of the metric is again a standard gluing argument.

\medskip

\epsfxsize =\textwidth
\epsfbox{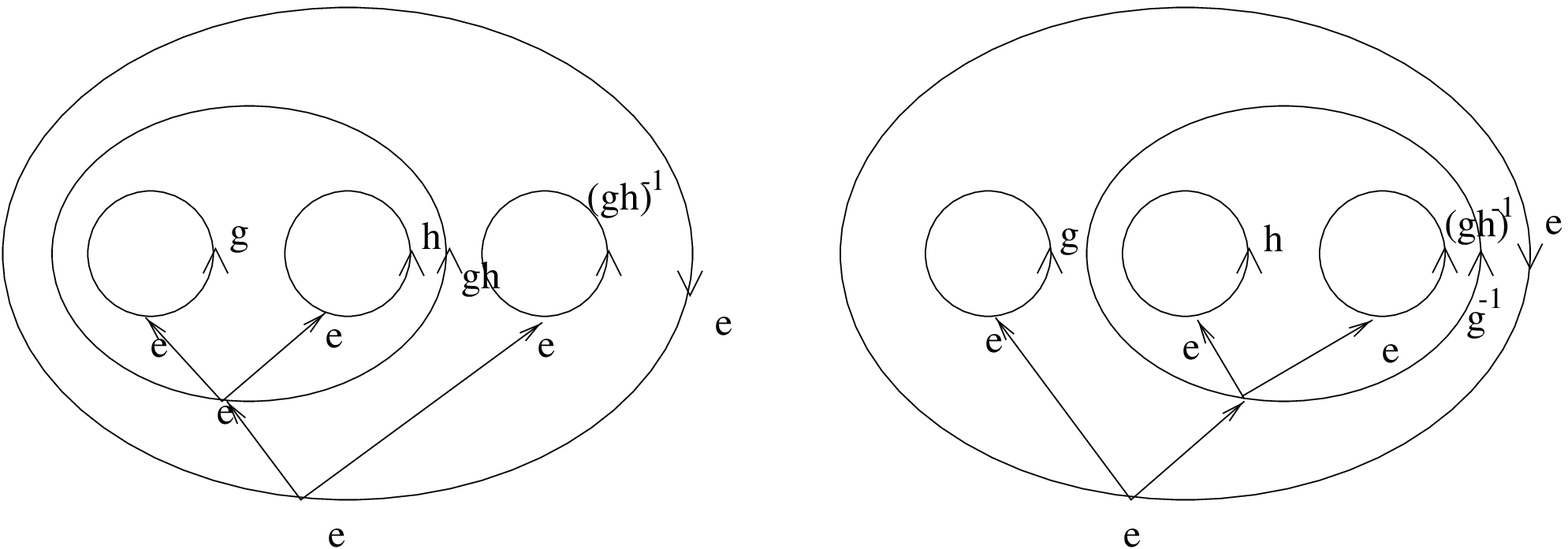}

\medskip

I.e.\ using $\bar\eps$ and associativity:

$$
\bar\eta(v_g\bar\circ v_h  ,v_{(gh)^{-1}})=
 \bar\eps((v_g\bar\circ v_h ) \bar \circ v_{(gh)^{-1}})=
\bar\eps(v_g\bar\circ (v_h \bar \circ v_{(gh)^{-1}}))=
\bar\eta(v_g, v_h \bar \circ v_{(gh)^{-1}})
$$

For axiom 1')

We use the following diagram:

$
\begin{array}
{@{\extracolsep{0pt}}c@{\extracolsep{0pt}}c@{\extracolsep{0pt}}
c@{\extracolsep{0pt}}c@{\extracolsep{0pt}}c@{\extracolsep{0pt}}}
(g,e)&\stackrel{id}{\rightarrow}&(g,e)&&(g,e)\\
\downarrow id&\Rightarrow&\downarrow \tilde C&\Rightarrow&\downarrow C\\
(g,e)&\stackrel{id}{\rightarrow}&(g,e)&&(g,g)\\
&&\Downarrow V&&\\
V_g&\stackrel{id}{\rightarrow}&V_g&\stackrel{id}{\rightarrow}&V_g\\
\downarrow id&&&&\downarrow \bar\varphi_g\\
V_g&\stackrel{id}{\rightarrow}&V_g&\stackrel{id}{\rightarrow}&V_g\\
\end{array}
$

where we used axiom viii) for the first move and axiom iii) for
the second.

\medskip

\epsfxsize =\textwidth
\epsfbox{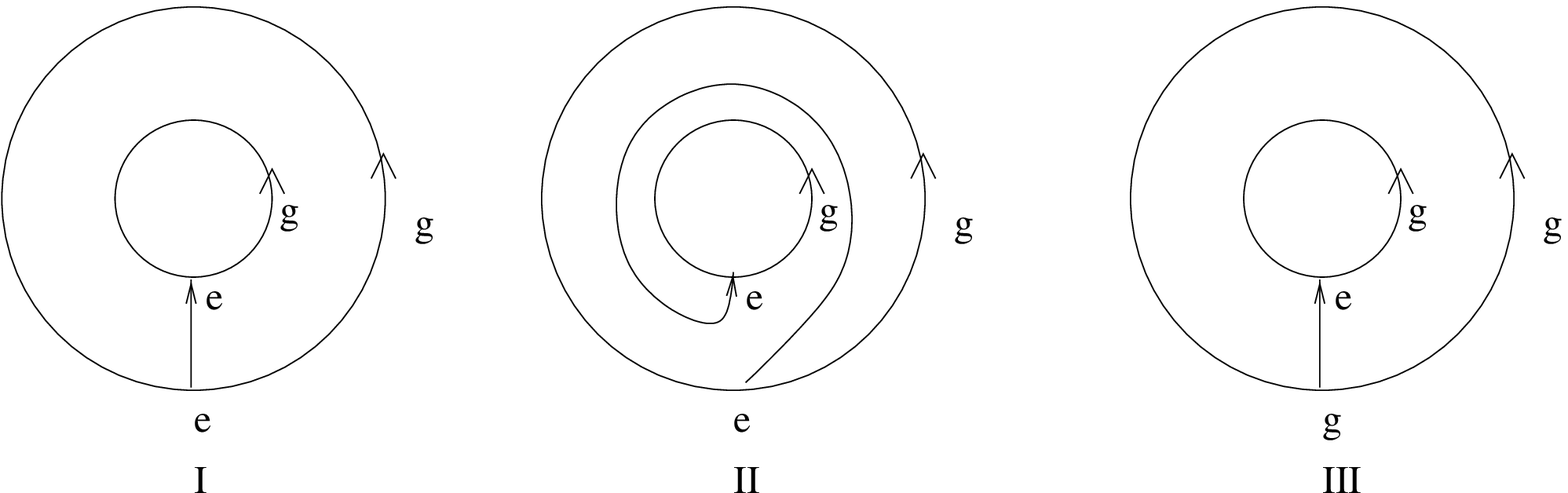}

\medskip

The axiom 2') follows from the diagrams

{\footnotesize
\begin{tabular}
{@{\extracolsep{0pt}}c@{\extracolsep{0pt}}c@{\extracolsep{0pt}}c@{\extracolsep{0pt}}c@{\extracolsep{0pt}}c@{\extracolsep{0pt}}c@{\extracolsep{0pt}}c@{\extracolsep{0pt}}}
$(g,e)\coprod(h,e)$&$\stackrel{(C_{(g,e)}^k, C_{(h,e)}^k)}{\longrightarrow}$&$(kgk^{-1},k)
\coprod (khk^{-1},k)$
&$\stackrel{T_{kgk^{-1},khk^{-1}}^{k}}{\longrightarrow}$&$(kghk^{-1},k)$&
$\stackrel{C_{(kghk^{-1},k)}^{k^{-1}}}{\longrightarrow}$&$(gh,e)$\\
$\downarrow$&&$\downarrow$&&$\downarrow$&&$\downarrow$\\
$V_g\otimes V_h$ &$ \stackrel {\bar\varphi_k\otimes\bar\varphi_k}{\longrightarrow}$&$V_{kgk^{-1}}\otimes V_{khk^{-1}}$&
$\stackrel{\circ^k}{\longrightarrow}$&$V_{kghk^{-1}}$&$\stackrel{\bar\varphi_{k^{-1}}}{\longrightarrow}$&$V_{gh}$\\
\end{tabular}
}

\medskip

\epsfxsize =4in
\epsfbox{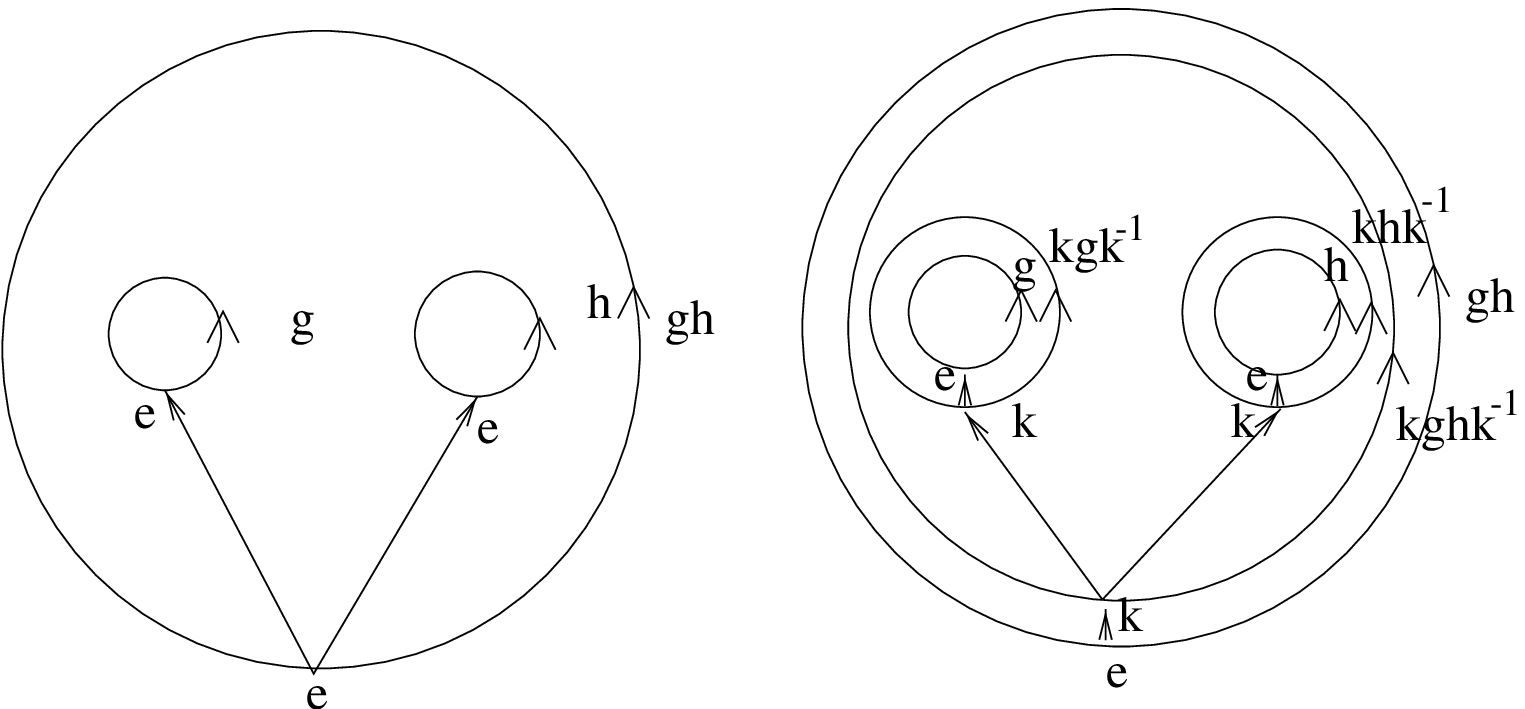}

\medskip

\begin{tabular}
{@{\extracolsep{0pt}}c@{\extracolsep{0pt}}c@{\extracolsep{0pt}}c@{\extracolsep{0pt}}c@{\extracolsep{0pt}}c@{\extracolsep{0pt}}c@{\extracolsep{0pt}}c@{\extracolsep{0pt}}}
$(g,e)\coprod(h,e)$&\multicolumn{5}{c}{$\stackrel{T^e_{gh}}{\longrightarrow}$}&$(gh,e)$\\
$\downarrow \|$&&&&&&$\downarrow \|$\\
$(g,e)\coprod(h,e)$&$\stackrel{(II_k,II_k)}{\longrightarrow}$&$(g,k)\coprod (h,k)$
&$\stackrel{T_{g,h}^{k}}{\longrightarrow}$&$(gh,k)$&
$\stackrel{II_{k^{-1}}}{\longrightarrow}$&$(gh,e)$\\
$\downarrow$&&$\downarrow$&&$\downarrow$&&$\downarrow$\\
$V_g\otimes V_h$ &$ \stackrel {\bar\chi_k\otimes \bar\chi_k}{\longrightarrow}$&$V_{g}\otimes V_{h}$&
$\stackrel{\circ^k}{\longrightarrow}$&$V_{gh}$&$\stackrel{\bar\chi_{k^{-1}}}{\longrightarrow}$&$V_{gh}$\\
$\downarrow \|$&&&&&&$\downarrow \|$\\
$V_g\otimes V_h$&\multicolumn{5}{c}{$\stackrel{\bar\circ}{\longrightarrow}$}&$V_{gh}$\\
\end{tabular}

The axiom 3') -- $G$--invariance of the metric -- follows from the following diagrams

\medskip

\epsfxsize = 4in
\epsfbox{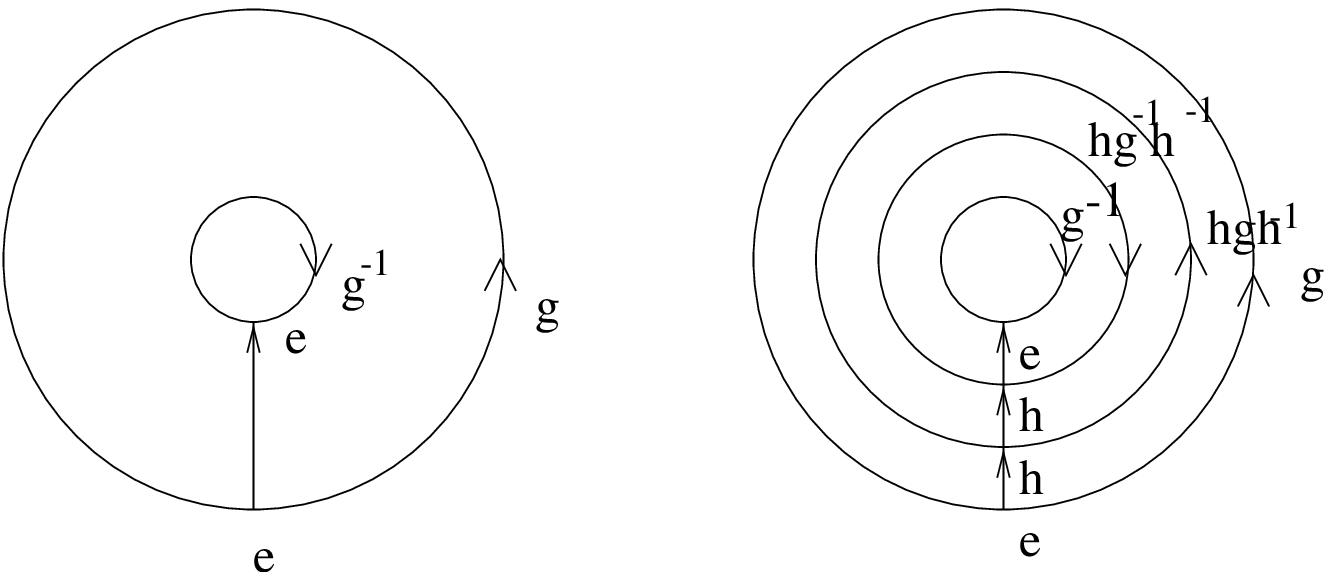}

\medskip

\begin{tabular}
{@{\extracolsep{0pt}}c@{\extracolsep{0pt}}c@{\extracolsep{0pt}}c@{\extracolsep{0pt}}c@{\extracolsep{0pt}}c@{\extracolsep{0pt}}c@{\extracolsep{0pt}}c@{\extracolsep{0pt}}}
$(g,e)\coprod\overline{(g,e)}$&\multicolumn{3}{c}{$\stackrel{A^e}{\longrightarrow}$}&$K$\\
$\downarrow \|$&&&&$\downarrow \|$\\
$(g,e)\coprod\overline{(g,e)}$&$\stackrel{(C,\bar C)}{\rightarrow}$&
$(kgk^{-1},k)\coprod \overline{ (kgk^{-1},k)}$
&$\stackrel{A^k}{\rightarrow}$&$K$\\
$\downarrow (id,\bar{} )$&&$\downarrow (id,\bar{} )$&&$\downarrow (id,\bar{} )$\\
$(g,e)\coprod(g^{-1},e)$&$\stackrel{(C, C)}{\rightarrow}$&
$(kgk^{-1},k)\coprod  (kg^{-1}k^{-1},k)$
&$\stackrel{A^k}{\rightarrow}$&$K$\\
$\downarrow$&&$\downarrow$&&$\downarrow $\\
$V_g\otimes V_{g^{-1}}$&$\stackrel{\bar\varphi_k \otimes \bar\varphi_k}{\rightarrow}$&
$V_{kgk^{-1}}\otimes V_{kg^{-1}k^{-1}} $
&$\stackrel{\eta^k}{\rightarrow}$&$K$
\end{tabular}

and

\begin{tabular}
{@{\extracolsep{0pt}}c@{\extracolsep{0pt}}c@{\extracolsep{0pt}}c@{\extracolsep{0pt}}c@{\extracolsep{0pt}}c@{\extracolsep{0pt}}c@{\extracolsep{0pt}}c@{\extracolsep{0pt}}}
$(g,e)\coprod\overline{(g,e)}$&\multicolumn{3}{c}{$\stackrel{A^e}{\longrightarrow}$}&$K$\\
$\downarrow \|$&&&&$\downarrow \|$\\
$(g,e)\coprod\overline{(g,e)}$&$\stackrel{(II_k,{II_k})}{\rightarrow}$&$(g,,k)\coprod \overline{ (g,k)}$
&$\stackrel{A^k}{\rightarrow}$&$K$\\
$\downarrow (id,\bar {})$&&$\downarrow (id,\bar{})$&&$\downarrow (id,\bar{})$\\
$(g,e)\coprod {(g^{-1},e)}$&$\stackrel{(II_k,{II_k}}{\rightarrow}$&$(g,,k)\coprod {(g^{_1},k)}$
&$\stackrel{A^k}{\rightarrow}$&$K$\\
$\downarrow$&&$\downarrow $&&$\downarrow$\\
$V_g\otimes V_{g^{-1}}$&$\stackrel{\bar\chi_k \otimes \bar\chi_{k^{-1}}}{\rightarrow}$&
$V_{g}\otimes V_{g^{-1}k} $
&$\stackrel{\eta^k}{\rightarrow}$&$K$\\
$\downarrow \|$&&&&$\downarrow \|$\\
$V_g\otimes V_{g^{-1}}$&\multicolumn{3}{c}{$\stackrel{\eta}{\longrightarrow}$}&$K$\\
\end{tabular}

Lastly axiom 4') comes from gluing a once punctured torus in two different ways.

\subsubsection{Proposition}
\label{reconstruction}
Given a $G$--Ramond algebra $V$  there is a unique
 $G_{\chi}$ orbifold theory $\cal V$ s.t.\ $V$ is its associated $G$--Ramond algebra.

{\bf Proof.}

We need to show that the data is sufficient to reconstruct the functor.
For the objects this is clear, due to the monoidal structure.
For discs, annuli and trinions the functor is defined by its basic ingredients, glueing
annuli $C_{g,h}$ to the basic trinion, and gauging with morphisms of type II.
For other morphisms of type I we notice that we can always decompose a surface into trinions, annuli
and discs. For each decomposition there are three choices.
The choice of a marking for the decomposition, a choice of orientation and a choice of pairs
of points over the marked curves. The second and the third choice can be seen to be irrelevant
by inserting two annuli $C_{e,k}$ and $A_{k,e}$ with suitable orientation in
a normal neighborhood of the curve in question.
The first choice is unique up to two operations [HT] which correspond to
associativity and the trace axiom, and is thus also irrelevant.

Combining the Propositions  \ref{fix} and \ref{reconstruction} with 
\ref{propdef} we obtain:
\subsubsection{Theorem}
\label{T1} There is a 1--1 correspondence between isomorphism
classes $G$-twisted Frobenius algebras and isomorphism classes of
$G_{\chi}$--orbifold theories as $\chi$ runs through the
characters of $G$.

{\bf Proof.} In the standard way, we make the $G$--twisted Frobenius algebras and the
 $G_{\chi}$--orbifold theories into categories by introducing the following morphisms.
For  $G$--twisted Frobenius algebras we use algebra homomorphisms respecting all the
additional structures and for  $G_{\chi}$--orbifold theories we use natural
transformations among
functors. The map of associating a $G$-twisted FA to a
 $G_{\chi}$--orbifold theory then turns into a full and faithful functor which is by
reconstruction surjective on the objects.

An immediate consequence of this is:
\subsubsection{Corollary} There is a 1--1 correspondence between  $G$-orbifold FA
and isomorphism classes of monoidal functors which are identity on
cylinders and satisfy the involutive property from $\mathcal
{GBCOB}$ to $\mathcal VECT$ which lift to
 $\mathcal {RGBCOB}$.

\subsection{Spectral flow}
In the previous paragraph, we chose the geometric version to
correspond to the Ramond picture, as is suggested by physics,
since we are considering the vacuum states in their
Hilbert--spaces at the punctures. From physics one expects that by
the spectral flow, the vacuum states should bijectively correspond
to the chiral algebra. In the current setting the difference only
manifests itself in a change of $G$--action given by a twist
resulting from the character $\chi$. We can directly produce this
$G$--action and thereby the $G$--Frobenius algebra by considering
the $G$--action not by $C_{e,k}^e$, but by $II_{k^{-1}} \circ
C_{e,k}^e$. This statement is proven by re-inspection of the
commutative diagrams in the proof of \ref{propdef}.

\section{Special $G$--Frobenius algebras}

In this section we restrict ourselves to a subclass of $G$--twisted Frobenius 
algebras. This subclass is large enough to contain all $G$--Frobenius algebras
 arising from singularities, symmetric products and spaces whose cohomology of
the fixed point sets are given by restriction of the comhomology of the ambient
space. The restriction allows us to characterize the possible $G$--Frobenius 
structures for a given collection of Frobenius algebras as underlying data in
terms of cohomological data. 
The restriction we will impose (cyclicity of the twisted sectors) can easily
be generalized to more generators; which will render everything matrix--valued.

\subsection{Definition} A {\it special $G$--Frobenius algebra}
is a $G$--twisted Frobenius algebra whose components $A_{g}$ are
cyclic $A_e$--modules, together with two collections of maps
$(r_{g}),(i_{g})$ indexed by $G$ where
\begin{itemize}
\item[] $r_{g}: A_e \rightarrow A_{g}$ is the  map of $A_e$--modules
induced by multiplication, we write  $a_{g} := r_{g}(a)$. In
particular, if one sets $1_{g} := r_{g}(1)$ we obtain $r_{g}(a) =
a 1_{g}$. Notice that it is equivalent to specifying the map
$r_{g}$ or generators $1_{g}$ of $A_{g}$ as cyclic $A_e$--modules.
\item[] $i_{g}: A_{g} \mapsto A_e$ are collections of maps  s.t.
each $i_{g}$ is an injection which splits the exact sequence
of $A_e$-modules.
\begin{equation}
0 \longrightarrow I_{g} \longrightarrow A
\raisebox {5pt}{\begin {array} [t]{c}
\stackrel {r_{g}} {\longrightarrow} \\
\stackrel{{\displaystyle \longleftarrow}} {{\scriptstyle  i_{g}}}
\end {array} }
A_{g} \longrightarrow 0
\end{equation}
\end{itemize}
Here $I_{g}\subset A_e$ is the annihilator of $1_{g}$ and thus of
$A_{g}$.
We denote the concatenation of $i_{g}$ and $r_{g}$ by $\pi_{g}$
$$
\pi_{g}=i_{g}\circ r_{g} : A_e \rightarrow A_e
$$
and we take the statement that $i_g$ is a section of $A_e$ modules to mean
$$
i_g(a_eb_g)= \pi_g(a_e)i_g(b_g)
$$

Furthermore the following should hold:

$$\forall g \in G: \varphi_{g}(1_{h}) = \varphi_{g,h}
1_{ghg^{-1}}$$ for some $\varphi_{g, h} \in K$ and
$\varphi_{g,e}=1$.

\subsubsection{Remarks}
Notice that
\begin{eqnarray*}
\pi_g(ab) &=& \pi_g(a)\pi_g(b)\\
r_g(ab) &=& \pi_g(a)r_g(b)
\end{eqnarray*}

In particular, multiplication by $1_{g}$ acts as projection:
$a 1_{g} = \pi_{g}(a) 1_{g}$
and $\pi_g(1)$ acts as identity on $A_g$.

\subsubsection{Special super $G$--Frobenius algebra}
The super version of special $G$--Frobenius algebras is straightforward. Notice that
since each $A_g$ is a cyclic $A_e$--algebra
its parity is fixed to be $(-1)^{\tilde g}:=\widetilde {1_g} $ times that of  $A_e$.
I.e.\ $a_g = i_{g}(a_{g})1_g$ and thus  $\widetilde{a_g} = \widetilde{i_g(a_{g})} \widetilde {1_g}$ .
In particular if $A_e$ is purely even, $A_g$ is purely of degree $\tilde g$.

\subsection{Proposition}
Let  $\eta_{g}:i_{g}(A_{g}) \otimes i_{g} (A_{g}) \rightarrow K$
be given by the following formula:
$$
\eta_{g}(a,b)  = \eta(r_{g}(a),r_{g^{-1}}(b))= \eps(\g_{g,g^{-1}}\pi_g(ab))
$$
Then  $i_{g}(A_{g})\subset A_e$
together with  $\eta_{g}$  and the induced multiplication $a \circ_g b
= i_g(a\circ r_g(b))=\pi_g(a\circ b)$
is a  Frobenius sub--algebra. Furthermore $\eta_{g}$ establishes
an isomorphism  between $i_{g}(A_{g})$ and
$i_{g^{-1}}(A_{g^{-1}})$.

{\bf Proof.}
That the multiplication is commutative and associative can be seen as follows

$$a \circ_g b= i_g(ab1_g)= i_g(ba1_g)=b\circ_g c$$

likewise

$(a \circ_g b ) \circ_g c= i_g((ab)c1_g)= i_g(a(bc)1_g)= a \circ_g (b  \circ_g c)$.
It is clear that
 $\eta_{g}$ is an invariant form.

The non--degeneracy follows directly from the non--degeneracy of
$\eta$: fix $a\in i_g(A_g)$ and let $a_{g} \in A_{g}$ and $b_{g^{-1}} \in A_{g^{-1}}$ s.t.
$i_{g}(a_{g})=a$ and
$\eta(a_{g},b_{g^{-1}}) \neq 0$. Set $b:= i_{g^{-1}}(b_{g^{-1}})$
and $\tilde b = \pi_{g}(b)$

Then
\begin{eqnarray*}
\eta_{g}(a,\tilde b)
&=& \eta(a1_{g},\tilde b1_{g^{-1}})
= \eta(\tilde b a1_{g},1_{g^{-1}})
= \eta (ba1_{g},1_{g^{-1}})\\
&=&  \eta( a 1_{g}, b1_{g^{-1}})
= \eta(a_{g},b_{g^{-1}})\neq 0
\end{eqnarray*}

The invariance follows from
$$
\eta_{g}(ab,c) = \eta(r_{g}(ab),r_{g^{-1}}(c))=
\eta(ab1_{g},c1_{g^{-1}}) = \eta(a1_{g},bc1_{g^{-1}})= \eta_{g}(a,bc)
$$

With these two properties, we get a map $i_{g}(A_{g})
\rightarrow i_{g}(A_{g})^*$ and $\eta$ gives an isomorphism
$i_{g}(A_{g})^* \rightarrow  i_{g^{-1}} (A_{g^{-1}})$
via $(a,b) \in A_{g} \otimes A_{g^{-1}} \mapsto \eta(a 1_{g},b
1_{g^{-1}})$.

\subsection{Proposition}
The following equality holds
$\pi_{g^{-1}}\circ\pi_{g}= \pi_{g}$
showing that $\pi_{g^{-1}}|_{i_g(A_g)}= id$
identifying $i_{g}(A_{g})$ and $i_{g^{-1}}(A_{g^{-1}})$.

\medskip

{\bf Proof.}
We have the following:
\begin{eqnarray*}
\eta(\pi_{g^{-1}}\circ\pi_{g}(a)1_{g},b_{g^{-1}})&=&
\eta(1_{g},
\pi_{g^{-1}}(\pi_{g}(a)i_{g^{-1}} (b_{g^{-1}}))1_{g^{-1}})\\
&=&
\eta(\pi_{g}(a)1_{g},i_{g^{-1}}(b_{g^{-1}})1_{g^{-1}})\\
&=& \eta(\pi_g(a)1_{g},b_{g^{-1}})
\end{eqnarray*}
which proves the statement due to the non--degeneracy of $\eta$.

\subsubsection{Remark}
We can also pull back $\eta_g$ and $\circ_g$ to $A_g$ which will make
$A_g$ into a Frobenius algebra and $i_g$ into an algebra homomorphism.

More precisely we put:
$$
a_g \circ_g b_g := r_g(i_g(a_g)i_g(b_g))= i_g(a_g)b_g
$$
and
$$
\eta_g(a_g,b_g):= \eps(i_g(a_g)i_g(b_g)\g_{g,g^-1})
$$

It then follows that
$$
i_g(a_g \circ_g b_g)= i_g(a_g) i_g(b_g)
$$

\subsection{Definition}  Let $A$ be a special $G$--twisted Frobenius
algebra
$A$. We define a {\it $G$ graded 2--cocycle with values in $A_e$} to be
a map $\g: G\times G \rightarrow A_e$ which satisfies
\begin{equation}
 \g_{g,h}:= \g(g,h) \in i_{gh}(A_{gh})
\end{equation}
and
\begin{equation}
\pi_{ghk}(\g_{g,h} \g_{gh,k})=
\pi_{ghk}(\g_{h,k}\g_{g,hk})
\end{equation}
We will call such a cocycle {\em associative} if also $\forall g,h,k \in G$:
\begin{multline}
\pi_{ghk}(\pi_{gh}(i_g(A_g)i_h(A_h))i_k(A_k)\g_{g,h}\g_{gh,k}) =\\
\pi_{ghk}(i_g(A_g)\pi_{hk}(i_h(A_h)i_k(A_k))\g_{h,k}\g_{g,hk})
\end{multline}
and we call a cocycle
{\em section independent} if $\forall g,h \in G$
\begin{equation}
(I_g+ I_h)\g_{g,h} \subset I_{gh}
\end{equation}

\subsection{Remark} Notice that if $g$ is section independent then
it follows that
the multiplication is independent of the choice of section and the cocycle is
automatically associative.

\subsection{Definition} A {\it non--abelian $G$ 2--cocycle with values in $k^{*}$} is a map
$\varphi: G\times G \rightarrow K$ which satisfies:
\begin{equation}
\label{grouphom}
\varphi_{gh,k}=\varphi_{g,hkh^{-1}} \varphi_{h,k}
\end{equation}
 where $\varphi_{g,h}:= \varphi(g,h)$ and
 $$\varphi_{e,g}=\varphi_{g,e}=1$$

Notice that in the case of a commutative group $G$ this says that
the $\varphi_{g,h}$ form a two cocycle with values in $K^*$.

Furthermore setting $g=h^{-1}$, we find:
$$
\varphi_{g^{-1},ghg^{-1}}=\varphi_{g,h}^{-1}
$$

\subsection{Proposition}  A special $G$--twisted Frobenius algebra $A$
gives rise to an associative $G$ 2--cocycle $\g$ with values in $A$ and
to a non--abelian $G$ 2--cocycle $\varphi$ with values in $k^{*}$,
which satisfy the compatibility equations
\begin{equation}
\label{grpcompat}
\varphi_{g,h}\g_{ghg^{-1},g} = \g_{g,h}
\end{equation}
and
\begin{equation}
\label{algaut}
\varphi_{k,g} \varphi_{k,h} \g_{kgk^{-1},khk^{-1}}
= \varphi_{k} (\g_{g,h}) \varphi_{k,gh}
\end{equation}

{\bf Proof.}

Given a special $G$ twisted Frobenius structure on $A$, we define
$\g_{g,h}\in i_{gh}(A_{gh})$ by
$$
1_{g}1_{h}= \g_{g,h} 1_{gh}
$$
and
$\varphi_{g,h} \in k^{*}$ by:
$$
\varphi_{g}(1_{h})= \varphi_{g,h}1_{ghg^{-1}}.
$$

By associativity of the multiplication, we find that the $\g_{g,h}$
define a graded 2--cocycle with values in $A$.
\begin{equation}
\label{ass}
\g_{g,h} \g_{gh,k}1_{ghk}
= \g_{g,h}1_{gh} 1_{k}
= (1_{g}1_{h}) 1_{k}
= 1_{g}(1_{h}1_{k})
= 1_{g}\g_{h,k}1_{hk}
= \g_{h,k}\g_{g,hk}1_{ghk}
\end{equation}
So that
$$
\pi_{ghk}(\g_{g,h} \g_{gh,k})= \pi_{ghk}(\g_{h,k}\g_{g,hk}).
$$
Redoing the calculation with general elements shows the associativity of the
cocycle.

$\varphi$ is a group homomorphism so that

$$
\varphi_{gh,k} 1_{ghkh^{-1}g^{-1}}=
\varphi_{gh}(1_{k}) = \varphi_{g}( \varphi_{h}(1_{k}))=
$$
$$
\varphi_{g}( \varphi_{h,k}1_{hkh^{-1}})
=\varphi_{g,hkh^{-1}} \varphi_{h,k} 1_{ghkh^{-1}g^{-1}}
$$
which yields
\begin{equation}
\varphi_{gh,k}=\varphi_{g,hkh^{-1}} \varphi_{h,k}
\end{equation}

\medskip

By $G$--twisted commutativity
$$
\g_{g,h}1_{gh} = 1_{g}1_{h}= \varphi_{g}(1_{h}) 1_{g}
= \varphi_{g,h}1_{ghg^{-1}} 1_g= \varphi_{g,h}\g_{ghg^{-1},g} 1_{gh}
$$

So $\g_{g,h}$ and the $\varphi_{g,h}$ satisfy
\begin{equation}
\varphi_{g,h}\pi_{gh}(\g_{ghg^{-1},g}) = \pi_{gh}\g_{g,h}
\end{equation}

$\varphi$ is also an algebra automorphism:

$$
\varphi_{k}(1_{g})\varphi_{k}(1_{h})=
\varphi_{k}(1_{g} 1_{h})
$$

Expressed in the $\varphi$'s and $\g$'s:
\begin{equation}
\varphi_{k,g} \varphi_{k,h} \g_{kgk^{-1},khk^{-1}}
= \varphi_{k} (\g_{g,h}) \varphi_{k,gh}
\end{equation}
which gives a formula for the action of $\varphi$ on the $\g$s.

\subsubsection{Corollary}
\label{list}
A special $G$-Frobenius algebra gives rise to collection of Frobenius--algebras
 $(A_g,\circ_g,1_g,\eta_g)_{g\in G}$ together with a $G$--action on $A_e$.
A graded $G$ 2--cocycle with values in $A_e$ and a
compatible non--abelian $G$ 2--cocycle with values in $k^{*}$.
Furthermore:
\begin{itemize}
\item[i)]
$\eta_{e}(\varphi_{g}(a),\varphi_{g}(b)) =
    \chi_{g}^{-2}\eta_{e}(a,b)$
\item[ii)] $\eta_{g}(a_{g},b_{g})=
\eta(i_{g}(a_{g})i_{g}(b_{g})\g_{g,g^{-1}},1)$
\item[iii)] The projective trace axiom
$\forall c \in A_{[g,h]}$ and $l_c$ left multiplication by $c$:
\begin{equation}
\chi_{h}\mathrm {Tr} (l_c  \varphi_{h}|_{A_{g}})=
\chi_{g^{-1}}\mathrm  {Tr}(  \varphi_{g^{-1}} l_c|_{A_{h}})
    \label{verlinde}
\end{equation}

\end{itemize}

\subsection{Proposition}(Reconstruction) Given a Frobenius algebra
$(A_e,\eta_{e})$ and a $G$-action on $A$ : $\varphi: G\times A \rightarrow A$
together with the following data
\begin{itemize}
\item [-] Frobenius algebras
$(A_{g},\eta_{g}|g \in G \setminus \{1\})$.
\item [-] Injective algebra homomorphisms $i_{g}: A_{g}\rightarrow A_e$
s.t.\ $i_{g}(A_{g})=i_{g^{-1}}(A_{g^{-1}})$

\item[-] Restriction maps: $r_g: A \rightarrow A_g$ s.t.\ $r_g \circ i_g =id$..
\item [-] A graded associative $G$ 2--cocycle $\gamma$ with values in $A_e$ and a
compatible non--abelian $G$ 2--cocycle  $\varphi$ with values in $k^{*}$, s.t.\
$\eta_g(a,b) = \eta_g(\g_{g,g^{-1}} i_g(a),i_g(b))$.
\item [-] A group homomorphism $\chi: G \rightarrow k^{*}$
\end{itemize}
such that  i)-iii) of \ref{list} hold,
then there is a unique extension of this set of data to a special $G$--twisted
Frobenius algebra, i.e. there is a unique special $G$--twisted
Frobenius algebra with these underlying data.

\medskip

{\bf Proof.}

Denote the unit of $A_e$ by $1$ and the unit of $A_{g}$ by $1_{g}$.

Denote by $\pi_{g}:A_e\rightarrow i_{g}(A_{g})$ the projection to
$A_{g}$: $\pi_g:= i_g \circ r_g$.

Let the $A_e$ module structure on $A_{g}$ be given by
$$
a b_{g}:= r_g(a i_g(b_g))
$$
It is clear that $A_{g}$ is a cyclic $A_e$--module generated by $1_{g}$.\\
The algebra structure is then determined by:
$$
a_{g}b_{h}:= i_{g}(a_{g})i_{h}(b_{h}) \g_{g,h}1_{gh}
$$
This multiplication is associative as can be seen by using the associativity of
the cocycle.

The $\varphi$'s determine the $G$-action by:
$\varphi_{g}(b_{h}) := \varphi_{g}(i_{h}(b_{h}))
\varphi_{g,h}1_{ghg^{-1}}$.

The compatibility and (\ref{algaut}) guarantee that this is a
representation and the $G$ action is
indeed an action by algebra automorphisms.

The $G$--twisted commutativity $a_{g}b_{h}= \varphi_{g}(b_{h}) a_{g}$
also follows from the compatibility.

The form $\eta$ is defined the following way:

$$
\eta(a_{g},h_{g^{-1}}):=
\eta_{e}(i_{g}(a_{g})i_{g^{-1}}(b_{a^{-1}})\g_{gg^{-1}},1)=
\eta_{g}(i_{g}(a_{g}),i_{g^{-1}}(b_{a^{-1}}))
$$
and
$$
\eta(a_{g},h_{h}):= 0 \; \mbox{ if }\; gh\neq 1
$$
This form is non--degenerate by assumption and it is invariant:

\begin{eqnarray*}
\eta(a_{g}b_{h},d_{(gh)^{-1}}) &=&
\eta_{e}(i_{g}(a_g)i_{h}(b_{h})i_{h^{-1}g^{-1}}(d_{h^{-1}g^{-1}})
\g_{gh}\g_{gh,h^{-1}g^{-1}},1)\\
&=&\eta_{e}(i_{g}(a_g)i_{h}(b_{h})i_{h^{-1}g^{-1}}
(d_{h^{-1}g^{-1}})\g_{h,h^{-1}g^{-1}}\g_{g,g^{-1}},1)\\
&=&\eta(g_{g},b_{h}d_{h^{-1}g^{-1}})
\end{eqnarray*}

\subsubsection{Remark}
By straightforward calculation we have that the projective trace axiom is equivalent to

$\forall g, h \in G; c 1_{ghg^{-1}h^{-1}} \in
A_{[g,h]}, c \in i_{hgg^{-1}h^{-1}}(A_{[gh]})$

\begin{multline}
\chi_h\varphi_{h,g}{\rm Tr}(l_{\g_{ghg^{-1}h^{-1},hgh^{-1}} c} \varphi_{h}|_ {i_{g}(A_{g})})=\\
\chi_{g^{-1}}\varphi_{g^{-1},ghg^{-1}}
{\rm Tr}(\varphi_{g^{-1}} l_{\g_{hg^{-1}h^{-1}g,g^{-1}hg})c} |_{ i_{h}(A_{h })})
\end{multline}

Notice that in the graded case (see below)
this condition only needs to be checked for $\deg(\g_{ghg^{-1}h^{-1},hgh^{-1}} c) \neq 0$
and $\deg(\g_{hg^{-1}h^{-1}g,g^{-1}hg}c) \neq 0$.

Furthermore if $[g,h]=e$ then
 $\g_{ghg^{-1}h^{-1},hgh^{-1}}=\g_{hg^{-1}h^{-1}g,g^{-1}hg}=1_e$ and
$\varphi_{g^{-1},ghg^{-1}}=\varphi_{g,h}^{-1}=\varphi_{g^{-1},h}$.

\subsection{Remarks}
\label{conds}
\begin{itemize}
\item[1)] If $A_{g} A_{h} \neq 0 $
the compatibility condition (\ref{grpcompat}) already
determines the $\varphi_{g,h} \in k^{*}$.
\item[2)] In particular: $\g_{g,g}= 0 $ unless $\chi_g=1$
and if $[g,h]=e$ we have $\varphi_{g,h}\varphi_{hg}=1$ or $\g_{g,h}=\g_{h,g}=0$.
\item[3)] If also
$A_{g}A_{h}A_{k}\neq 0$
the elements defined by  (\ref{grpcompat})  automatically satisfy the conditions of non--abelian 2--cocycles
and the condition (\ref{algaut}) is automatically satisfied.
\end{itemize}

\medskip

{\bf Proof}
Without loss of generality, we may assume
that $1_{g}1_{h}1_{k} \neq 0$ then due to the condition (\ref{grpcompat})
$$
1_{g}1_{h}1_{k} = (1_{g} 1_{h}) 1_{k}
= \varphi_{gh}(1_{k}) (1_{g} 1_{h})
=(\varphi_{gh}(1_{k}) 1_{g}) 1_{h}
=\varphi_{gh,k} 1_{ghkh^{-1}g^{-1}}1_{g} 1_{h}
$$
and using associativity, we similarly obtain
$$
1_{g}1_{h}1_{k} = 1_{g} (1_{h} 1_{k})
=1_{g} \varphi_{h}(1_{k}) 1_{h}
= \varphi_{g}(\varphi_{b,k}1_{hkh^{-1}}) 1_{g} 1_{h}
= \varphi_{g,hkh^{-1}}\varphi_{h,k} 1_{ghkh^{-1}g^{-1}}1_{g} 1_{h}
$$

For the first statement in 2) one just needs to plug $g=h$ into (\ref{grpcompat})
and for the second one has to apply the formula twice.
For 3)  notice that
$1_{k}1_{g}1_{h} = 1_{k}(1_{g}1_{h})= \varphi_{k}(1_{g}1_{h}) 1_{k}$
and on the other hand
$1_{k}1_{g}1_{h} = \varphi_{k}(1_{g})1_{k}1_{h})
=\varphi_{k}(1_{g})\varphi_{k}(1_{h})1_{k}$.

A useful technical Lemma to show that $\g_{g,h}\neq 0$ is the following

\subsection{Lemma}
\label{zerocheck}
If $\g_{g,h}=0$ then $\pi_{h}(\g_{g,g^{-1}})=0$ and
$\pi_{g}(\g_{h,h^{-1}})=0$

{\bf Proof.}

If $\g_{g,h}=0$ then\\
$0=\pi_{h}(\g_{g^{-1},gh}\g_{g,h}))=
\pi_{h}(\g_{{g^{-1}},g}\g_{e,h}) = \pi_{h}(\g_{g^{-1},g}) =
\pi_{h}(\g_{g,g^{-1}})$

and also

$0=\pi_{g}(\g_{g,h}\g_{gh,h^{-1}}))=
\pi_{g}(\g_{g,e}\g_{h^{-1},h}) = \pi_{g}(\g_{h,h}^{-1}) $

\subsubsection{Graded special $G$--Frobenius algebras}

Consider a set of graded Frobenius algebras satisfying the reconstruction data:
$\{(A_{g},\eta_{g}):g \in G\}$  with degrees $d_{g}:=\deg(\eta_{g})$\\
s.t. $A_{g}\simeq A_{g^{-1}}$.  E.g.\ in the cohomology of fixed
point sets $d_{g}$ is given by the dimension and for the Jacobian
Frobenius manifolds (see the next section) $d_{g}$ fixed by the
degree of $Hess(f_{g}) = \rho_{g}$. Furthermore, the reconstructed
$\{\eta|_{(A_{g}\otimes A_{g^{-1}}}, g \in G\}$ have degree
$d_{g}=d_{g^{-1}}$.

For a $G$--twisted FA the degrees all need to be equal to $d:=d_{e}$.
To achieve this, one can shift the grading in each
$A_{g}$ by $s_{g}$. This amounts to assigning degree $s_{g}$ to $1_{g}$.
This is the only freedom, since the multiplication should be degree--preserving and
all $A_g$ are cyclic.

Set\\
$s_{g}^+ := s_{g}+s_{g^{-1}}$\\
$s_{g}^-:= s_{g}-s_{g^{-1}}$

Then $s_{g}^+:= d-d_{g}$ for grading reasons, but the shift $s^-$ is more elusive.
\subsubsection{Definition}
The standard shift for a Jacobian Frobenius algebra
is given by
$$s_{g}^+:= d-d_{g}$$
and
\begin{multline*}
s_{g}^- := \frac{1}{2\pi i}\mathrm{tr} (\log(g))-\mathrm{tr}(\log(g^{-1})):=
\frac{1}{2\pi i}(\sum_i \l_i(g)-\sum_i \l_i(g^{-1}))\\
=\sum_{i: \l_i \neq 0} 2 (\frac{1}{2\pi i}\l_i(g)-1)
\end{multline*}
where the $\l_i(g)$ are the logarithms of the eigenvalues of $g$ using the branch with values in $[0,2\pi)$ i.e.\
cut along the positive real axis.

In this case we obtain:

$$
s_{g}= \frac{1}{2}(s_g^+ + s_g^-)= \frac{1}{2}(d-d_g) + \sum_i (\frac{1}{2\pi i}\l_i(g)-\frac{1}{2})
$$

\subsubsection{Remark}
The shift $s_{g}^-$ is canonical in the case of quasi--homogeneous
singularities upon replacing the classical monodromy operator $J$ by
$Jg$. This will be discussed elsewhere [K2].

The degree of $\g_{g,g^{_1}}$ is $s^+$ from comparing degrees in the
equation $1_g 1_{g^{-1}}= \g_{g,g^{_1}}1_e$.

\subsubsection{Reconstruction for graded special $G$--Frobenius algebras}
In the Reconstruction program the presence of a non--trivial grading can greatly
simplify the check of the trace axiom. E.g.\ if $A_{[g,h]}$ has no element of degree $0$ then
both sides of this requirement are $0$ and if $[g,h]=e$ one needs only to look at the
special choices of
$c$ with $\deg(c)=0$ which most often is just $c=1$, the identity.

\subsubsection{Ramond--grading}
The grading in the Ramond--sector is by the following definition

$$\deg(v):=-\frac{d}{2}$$

This yields

$$\deg(\bar\eta) = 0 \mbox{ and } \deg(\bar\circ)=\frac{d}{2}$$

\section{Jacobian Frobenius algebras}

\subsection{Definition} A Frobenius algebra $A$ is called {\it Jacobian}

if it can be represented as the Milnor ring of a function $f$. I.e. if
there is a function $f\in {\Cal O}_{{\bf A}^{n}_{K}}$
s.t. $A = {\Cal O}_{{\bf A}^{n}_{K}}/J_{f}$
 where $J_{f}$ is the Jacobian
ideal of $f$. And the bilinear form is given by the residue
pairing. This is the form given by the the Hessian of $f$:
$\rho={\rm Hess}_{f}$.

If we write ${\Cal O}_{{\bf A}^{n}_{k}}= k[x_{1}\dots x_{n}]$,
$J_{f}$  is  the ideal spanned by the $\frac{\del f}{\del x_{i}}$.

A {\it realization of a Jacobian Frobenius algebra} is a pair
$(A,f)$ of a Jacobian Frobenius algebra and a function $f$  on
some affine $k$ space ${\bf A}_{k}^{n}$, i.e. $f \in {\Cal
O}_{{\bf A}_{k}^{n}}= k[x_{1}\dots x_{n}]$ s.t. $A= k[x_{1}\dots
x_{n}]$ and $\rho := {\rm det}(\frac{\del^{2} f}{\del x_{i}\del
x_{j}})$.

A {\it small realization of a Jacobian Frobenius algebra} is a
realization of minimal dimension, i.e. of minimal $n$.

\subsection{Definition} A {\it natural $G$ action on a realization of a
Jacobian Frobenius algebra $(A_e,f)$ } is a linear $G$ action on
${\bf A}_{k}^{n}$ which  leaves $f$ invariant.

Given a natural $G$ action on a realization of a Jacobian
Frobenius algebra $(A,f)$ set for each $g\in G$, ${\Cal O}_{g}:=
{{\Cal O}_{ {\rm Fix}_{g}({\bf A}_{k}^{n})}}$. This is the ring of
functions of the fixed point set of $g$ for the $G$ action on
${\bf A}_{k}^{n}$. These are the functions fixed by $g$: ${\Cal
O}_{g}= k[x_{1},\dots, x_{n}]^{g}$.

Denote by $J_{g}:= J_{f|_{{\rm Fix}_{g}({\bf A}_{k}^{n})}}$ the
Jacobian ideal of $f$ restricted to the fixed point set of $g$.

Define
\begin{equation}
A_{g}:= {\Cal O}_{g}/J_{g}
\end{equation}
The $A_{g}$ will be called twisted sectors for $g \neq 1$.
Notice that each $A_{g}$ is a Jacobian Frobenius algebra with the
natural realization given by $(A_{g}, f|_{\mathrm {Fix}_{g}})$.
In particular, it comes equipped with an invariant bilinear form
$\tilde\eta_{g}$
defined by the element $\mathrm {Hess}(f|_{\mathrm {Fix}_{g}})$.

For $g = 1$ the definition of $A_e$ is just the
realization of the original Frobenius algebra,
which we also call the untwisted sector.

Notice there is a restriction morphism $r_{g}: A_e \rightarrow
A_{g}$ given by 
$a \mapsto a|_{\mathrm{Fix}_g} \;{\rm mod} \;J_{g}$.

Denote $r_{g}(1)$ by $1_{g}$. This is a non--zero element of
$A_{g}$ since the action was linear.
Furthermore it generates $A_{g}$ as a cyclic $A_e$ module.

The set ${\rm Fix}_{g}{\bf A}_{k}^{n}$ is a linear subspace. Let
$I_{g}$ be the vanishing ideal of this space.

We obtain a sequence
$$
0 \rightarrow I_{g} \rightarrow A_e\stackrel {r_{g}}{\rightarrow}
A_{g}\rightarrow 0
$$

Let $i_{a}$ be any splitting of this sequence induced by the inclusion:
$\hat i_{g}: {\Cal O}_{g} \rightarrow {\Cal O}_{e}$ which descends due
to the invariance of $f$.

In coordinates, we have the following description. Let ${\rm
Fix}_{g}{\bf A}_{k}^{n}$ be given by equations $x_{i}=0: i \in
N_{g}$ for some index set $N_{g}$.

Choosing complementary generators $x_{j}: j \in T_{g}$ we have
${\Cal O}_{g}= k[x_{j}:j \in T_{g}]$ and ${\Cal O}_{e}=
k[x_{j},x_{i}:j \in T_{g}, i\in N_{g}]$. Then $I_{g}=(x_{i}: i \in
N_{g})_{{\Cal O}_{e}}$ the  ideal in ${\Cal O}_{e}$ generated by
the $x_{i}$ and ${\Cal O}_{e} = I_{g}\oplus  i_{g}(A_{g})$ using
the splitting $i_{g}$ coming from the natural inclusion $\hat
i_{g}:k[x_{j}:j \in T_{g}] \rightarrow k[x_{j},x_{i}:j \in T_{g},
i\in N_{g}]$. We also define the projections
$$
\pi_{g}: A_1 \rightarrow A_{g}; \pi_{g} = i_{g} \circ r_{g}
$$
which in coordinates are given by $f \mapsto f|_{x_{j}=0: j \in
N_{g}}$
Let
$$
A := \bigoplus_{g \in G} A_{g}
$$

where the sum is a sum of $A_e$ modules.

\subsection{Definition}
A {\em discrete torsion} for a group $G$ is a map from commuting
pairs $(g,h) \in G \times G: [g,h]=e$ to $k^*$ with the
properties:

\begin{equation}
\eps(g,h)=\eps(h^{-1},g) \quad \eps(g,g)=1 \quad \eps(g_1g_2,h)=
\eps(g_1,h)\eps(g_2,h)
\end{equation}

\subsection{Theorem}
Given a natural $G$ action on a realization of a Jacobian
Frobenius algebra $(A_e,f)$ with a quasi--homogeneous function $f$
of degree $d$ and type ${\bf q}= (q_{1}, \ldots , q_{n})$ together
with the natural choice of splittings $i_{g}$ the possible
structures of special $G$ twisted Frobenius algebra on the $A_e$
module $A := \bigoplus_{g \in G} A_{g}$ are in 1--1 correspondence
with elements of $\bar Z^{2}(G, A)$, where $\bar Z^{2}$ are $G$
graded cocycles and a compatible non--abelian two cocycle
$\varphi$ with values in $k^{*}$, which define a choice of
discrete torsion.
\medskip

{\bf Proof.}
By the reconstruction Theorem the structure of a special $G$ twisted
Frobenius algebra on the given data $A_{g}, \eta_{g},i_{g}$
is determined by:

\begin{itemize}
\item [-] A graded $G$ 2--cocycle with values in $A_e$ and a
compatible non--abelian $G$ 2--cocycle with values in $k^{*}$
\item [-] A group homomorphism $\chi: G \rightarrow k^{*}$
\end{itemize}

subject to the conditions:

\begin{itemize}
\item[i)]
$
\eta_{e}(\varphi_{g}(a),\varphi_{g}(b)) =
    \chi_{g}^{-2}\eta_{e}(a,b)
$
\item[ii)] $\eta_{g}(a_{g},b_{g})=
\eta(i_{g}(a_{g}b_{g})\g_{gg^{-1}},1)$
\item[iii)] For any $c \in A_{[gh]}$:
\begin{equation*}
\chi_{h}{\rm STr} (l_c \varphi_{h}|_{A_{g}})
= \chi_{g^{-1}}{\rm STr} (\varphi_{g^{-1}} l_c|_{A_{h}} )
\end{equation*}
where we can restrict to the $c$ with $\deg(c)=0$.
\end{itemize}

The $\chi_{g}$ are fixed by the equation

$$
(-1)^{\tilde g}\dim(A_g)= {\rm STr} (id|_{A_{g^{-1}}})
= \chi_{g}{\rm STr} (\varphi_{g} |_{A_e})=\chi_{g}{\rm Tr} (\varphi_{g} |_{A_e})
$$
To calculate the trace on the RHS, we use the character function for a
morphism $g$ of degree $0$ on graded module $V = \oplus_{n}V_{n}$:

$$
\chi_{V_{n}}(g,z):= \sum_{n,\mu} \mu \dim (V_{\mu,n}) z^{n}
$$
where $V_{\mu,n}$ is the Eigenspace of Eigenvalue $\mu$ on the
space $V_{n}$. We will use the grading induced by the
quasi--homogeneity. I.e. let $N$ be such that $q_{i}= Q_{i}/N$
with $Q_{i}\in {\bf N}$ and $N$ s.t. $|G| {\big |} N $. Then a
monomial has degree $n$ if its quasi homogeneous degree is $n/N$.
This is the natural grading for the quasi--homogeneous map
$\mathrm {grad} (f)$. Notice that since $g$ commutes with $f$ it
preserves the grading. It is clear that this character behaves
multiplicatively under concatenations of quasi--homogeneous
functions. Therefore by applying Arnold's method, we can pass to a
cover of $k^{n}$ with the projection map $T:
T(x_{1},\ldots,x_{n})= (x^{q_{1}},\ldots,x^{q_{n}})$ and calculate
the character for $T$ and for ${\rm grad} (f) \circ T$. Then
repeating the argument in a simultaneous Eigenbasis of $g$ and the
grading of [A], we obtain:

$$
\chi_{A_e}(g,z)= \prod_{i=1}^{n}
\frac{(\tilde\mu_{i}z)^{N-Q_{i}}-1}{(\tilde \mu_{i}z)^{Q_{i}}-1}
$$

where the $\tilde \mu_{i}$ are the Eigenvalues of some lift of
the action of $g$ i.e. $\tilde\mu_{i}^{Q_{i}}=\mu_{i}$.
Notice that since $|G| {\big |}N$ $\tilde\mu_i^{N}=1$, so that in the limit
of $z \rightarrow 1$, we obtain:

$$
{\rm Tr}(\varphi_{g} |_{A_e})
 =  \prod_{i:\mu_{i}\neq 1} -\mu_{i}^{-1}\prod_{i:\mu_{i}=1}\frac{1}{q_{i}}-1
= (-1)^{|N_{g}|} \mathrm{det}(g)^{-1} \;\mathrm{dim} (A_{g})
$$
so that
$$
\chi_{g}= (-1)^{\tilde g} (-1)^{|N_{g}|} \mathrm{det}(g)
$$
We set
$$
\sigma(g):= \tilde g +|N_g| \mod 2
$$
and call it the sign of $g$.
Notice that $\sigma \in \mathrm {Hom}(G,\Z2)$ since both $\det$ and $\chi$ are characters.
Also notice that a choice of sign corresponds to a choice of parity and vice--versa.
Thus we obtain
$$
\chi_{g}= (-1)^{\sigma(g)} \det(g)
$$

This ensures condition i). Since  for
$\rho = \mathrm{Hess}(f)$:
$$
 \varphi_{g}(\rho) = \mathrm {det}(g)^{-2}\rho
$$

We  rescale $\eta_{g}$ by
$((-1)^{\tilde g}\chi_g)^{1/2}\eta_g$
where we choose to cut the plane along the negative real axis.
This uniquely defines a square root unless $(-1)^{\tilde g}\chi_g=-1$.
In the case that
$g^2\neq e$ we can choose roots $i$ and $-i$ for $g$ respectively
 $g^{-1}$. The only case that has no solution would be the case of $g^2 =e$ and
$(-1)^{\tilde g}\chi_g=-1$, but this means that
either $\chi_g =-1$ and $\tilde g =1$ or
$\chi_g =-1$ and $\tilde g =0$ which cannot happen, since in this case
$(-1)^{|N_g|} =\det(g)$ and $(-1)^{\tilde g}(-1)^{|N_g|}\det(g)=\chi_g$.
Then $\eta_g$ will satisfy the reconstruction conditions.

Finally, we need to check the validity of iii). Notice that since the
multiplication is graded the traces are $0$ unless ${\rm deg}(\g_c)=0$
so that we can assume that $\g_c= 1$. In this case, we have to show:
$$\chi_{h}{\rm STr} (\varphi_{h}|_{A_{g}})
= \chi_{g^{-1}}{\rm STr} (\varphi_{g^{-1}}|_{A_{h}})$$

Let $x_{i}$ be a basis of $A_e$ in which $g$ is diagonal. Then we
have to compute the trace of the action of $h$ on the sub-algebra
generated by the $x_{i}$ with eigenvalue 1 under the action of
$g$. This is just the truncated version of the calculation above,
so diagonalizing $h$ on $k[x_{i}:i \in T_{g}]$ we find using the
same characteristic functions and rationale as before:

\begin{eqnarray*}
&&\chi_h {\rm STr} (\varphi_{h}|_{A_{g}})=
\chi_h \varphi_{h,g}(-1)^{\tilde g}{\rm Tr} (\varphi_{h}|_{i_{g}(A_{g}) })\\
& =& \chi_h \varphi_{h,g}(-1)^{\tilde g}\prod_{j: \nu_{j}\neq 1} -\nu_{j}^{-1}
\prod_{j:\nu_{j}=1}\frac{1}{q_{j}}-1\\
&=& \chi_h \varphi_{h,g}(-1)^{\tilde g}(-1)^{|T_{g}|}(-1)^{|T_{g}\cap T_{h}|  } \det (h|_{T_{g}})^{-1}
\dim(i_g(A_{g})\cap i_h(A_{h})) \\
&=&  \varphi_{h,g}(-1)^{\sigma(g)} (-1)^{\sigma(h)}  \mathrm{det}(h)
\det (h|_{T_{g}})^{-1}
(-1)^{N}(-1)^{|T_{g}\cap T_{h}|  } \dim(i_g(A_{g})\cap i_h(A_{h}))\\
&=& \varphi_{h,g} (-1)^{\sigma (hg)}\det (h|_{N_{g}})
 (-1)^{|T_{g}\cap T_{h}| +N} \dim(i_g(A_{g})\cap i_h(A_{h})) \\
&=&\eps(h,g)  T(h,g)
\end{eqnarray*}
where $\nu_{j}$ are the Eigenvalues of $h$ on $i_{g}(A_{g})$
and

$$\eps(g,h)= \varphi_{g,h} (-1)^{\sigma(g)\sigma(h)}\det (g|_{N_{h}})
$$
and we set $\det(g|_{N_h}) := \det(g)\det^{-1}(g|_{T_{h}})$
if $[g,h]\neq e$.
$$
 T(h,g)= (-1)^{\sigma(g)\sigma(h)} (-1)^{\sigma(g)+\sigma(h)} (-1)^{|T_{g}\cap T_{h}|  +N} \dim(i_g(A_{g})\cap i_h(A_{h}))
$$
we have that
$$
T(h,g)=T(g,h)=T(g^{-1},h)
$$

Notice also that if $[g,h]=e$
$$
\eps(gh,k) = \eps(g,k)\eps(h,k)
$$

and

$$
T(gh,h) = T(hg,h)=T(g,h)
$$

On the other hand

\begin{eqnarray*}
&&{\rm STr} (\varphi_{g^{-1}}|_{A_{h}})=
\eps(g^{-1},h)(-1)^{\sigma(g^{-1}h)} (-1)^{|T_{g}\cap T_{h}| +N }
\dim(i_g(A_{g})\cap i_h(A_{h})) \\
&=&\varphi_{h,g} (-1)^{\sigma (g^{-1}h)}\det (g^{-1}|_{N_{h}})
(-1)^{\tilde g \tilde h} (-1)^{|T_{g}\cap T_{h}| } \dim(i_g(A_{g})\cap i_h(A_{h})) \\
&=&\eps(g^{-1},h) T(g^{-1},h)
\end{eqnarray*}
where $\mu_{j}$ are the Eigenvalues of $g^{-1}$ on $i_{h}(A_{h})$.

Finally we see that the $\varphi$ are determined by the $\eps(g,h)$ which have to satisfy

$$
\eps(g,h)=\eps(h^{-1},g) \quad \eps(g,g)=1 \quad \eps(g_1g_2,h)= \eps(g_1,h)\eps(g_2,h)
$$
which means that the choices of the projective factors for
the $G$--action correspond in a 1--1 fashion to a choice of
discrete torsion.

\subsubsection{Remarks}
\begin{itemize}
\item[1)] For a section independent cocycle $\g$
we have that
$$
\eta(\g_{g,g^{-1}}, I_g)= 0 \text { and } \eta(\g_{g,g^{-1}},i_g(\rho_g)) =1
$$
where $\rho_g$ is the defining element of $\eta_g$.
Thus $\g_{g,g^{-1}}$ is uniquely defined as
$\g_{g,g^{-1}}= \check{r}_{g}(1_{g})$.

where $\check{}$ is the dual linear map in the sense of vector-spaces
with non-degenerate bilinear forms.

Furthermore if $\g_{g,h}=\g_{h^{-1},g^{-1}}=1$ then
$$
\g_{gh,(gh)^{-1}} = \g_{g,g^{-1}} \g_{h,h^{-1}}
$$

\item[2)]
If $\forall g,h, \in G:J_h = J_g$
any choice of cocycle $\g$ with values in $k^{*}$ will give a special
G--twisted FA. There is only one choice of compatible non--abelian
cocycle $\varphi$.

\item[3)] For any other solution $\g$ by (\ref{conds}), we must have
$\g_{g,h}=0$ if $[g,h]=e$  unless $\varphi_{g,h}\varphi_{h,g}=1$
or in other words
$$
\det(g|_{N_h})\det(h)|_{N_g}=1
$$
since $\eps(g,h)=\eps(h^{-1},g)=\eps(h,g)^{-1}$
\end{itemize}

{\bf Proof.}
For 1) notice that
$$
\eta(\g_{g,g^{-1}}, I_g)= \eta(\g_{g,g^{-1}} I_g,1)=0
$$
and
$$
\eta(\g_{g,g^{-1}},i_g(\rho_g)) =\eta_g(1,\rho_g)=1
$$
regardless of the choice of $i_g$, by the previous equation.
For the uniqueness notice that in the graded case
the degree $d-\deg(\g_{g,g^{-1}})$ piece of $A$ is the direct sum of $i_g(\rho_g)$ and
the degree $d-\deg(\g_{g,g^{-1}})$ piece of $I_g$.

Furthermore we have that
$$
\g_{gh,(gh)^{-1}} = \g_{g,h}\g_{gh,(gh)^{-1}}
= \g_{g,g^{-1}}\g_{h,h^{-1}g^{-1}}
= \g_{g,g^{-1}} \g_{h,h^{-1}}
$$
since
$$
\pi_{g^{-1}} (\g_{h,h^{-1}g^{-1}}) = \pi_{g^{-1}} (\g_{h^{-1},g^{-1}} g_{h,h^{-1}g^{-1}})
=\pi_{g^{-1}}( \g_{h,h^-1}\g_{e,g^{-1}})
$$
and $I_{g^{-1}}\g_{g,g^{-1}}=0$. This is useful for calculations; another proof is: $\forall a\in
A_e$
$$
\eta(a,\check r_g(1_g)) = \eta_g(r_g(a),1_g)= \eta(\pi_g(a),\g_{g,g^{-1}})=\eta(a,\g_{g,g^{-1}})
$$

For 2) notice that all $\pi_{h}(\g_{g,g^{-1}})\neq 0$. Thus by
\ref{zerocheck} the $\g_{g,h}\neq 0$ and furthermore since $r_{g}=id$,
we see that the $\g_{g,h}\in k^{*}$.

For 3) we use the compatibility twice. If $[g,h]=e$:\\
$\g_{g,h}=\varphi_{g,h}\g_{h,g}=\varphi_{g,h}\varphi_{h,g}\g_{g,h}$

Now since by assumption $\g_{g,h}\neq 0$ we get the desired result.

\section{Mirror construction for special $G$--Frobenius algebras}

\subsection{Double grading}
We consider Frobenius algebras with grading in some abelian group $I$.
$$A= \bigoplus_{i\in I} A_i$$

This grading can be trivially extended to a double grading with values in $I\times I$ in two ways

$$
A^{cc} =  \bigoplus_{i\in I} A_{i,i}
$$
and
$$
A^{ac} =  \bigoplus_{i\in I} A_{i,-i}
$$

corresponding to the diagonal $\Delta: I \rightarrow I\times I$ and $(id,-)\circ \Delta: I \rightarrow I\times I$.
We call bi--graded Frobenius algebras of this form of $(c,c)$--type and of $(a,c)$--type, respectively.
In the language of Euler fields we consider the field $(E, \bar E)= (E,E)$ or  $(E, \bar E)= (E,-E)$.

These gradings become interesting for special $G$--Frobenius algebras, since in that case the shifts
will produce a possible non--diagonal grading.

\subsection{Definition}
Given a graded special $G$--Frobenius algebra we assign the following bi--degrees to $1_g$

$$
(E,\bar E) (1_g) := (s_g,s_g^{-1})
$$

It is clear that $A_e$ is of $(c,c)$--type. $A$ is however only of $(c,c)$ type if $s_g = s_g^{-1}$.

Furthermore for the Ramond--space of $A$ we assign the following bi--degree to $v$

$$
(E,\bar E) (v):= (-\frac{d}{2},-\frac{d}{2})
$$

\subsection{Euler--twist (spectral flow)}

In this section, we consider a graded special $G$--Frobenius algebra and construct a new
vector--space from it. We denote the grading operator by $E$:

\subsubsection{Definition}
The twist--operator $j$ for an Euler--field $E$ is

$$j:=\ exp(2\pi i E)$$

We denote the group generated by $j$ by $J$.

We call a special $G$--Frobenius {\em Euler} if there is a
special $\tilde G$--Frobenius algebra $\tilde A$ of which $A$ is a
subalgebra where $\tilde G$ is a group that has $G$ and $J$ as
subgroups.

\subsubsection{Definition}
The dual  $\check A$ to an Euler  special $G$--Frobenius algebra
$A$ of $(c,c)$--type  is the vector space
$$
\check A := \bigoplus_{g\in G} \check A_g := \bigoplus_{g\in G}
V_{gj^{-1}}
$$
with the $A_e$ and $G$--module structure determined by  $\tau_j^g:
\check A_g \simeq V_{gj^{-1}}$ together with the bi--grading
$$
(E,\bar E) (\check1_g):= (s_{gj^{-1}}-d,s_{gj^{_1}})
$$
where $\check 1_g$ denote the generator of $\check A_g$ as
$A_e$--module and the bi--linear form
$$
\check\eta := \tau_j^*(\bar\eta)
$$
where $\tau_j := \oplus_{g\in G}\tau_j^g$ and $V$ and $\bar \eta$
refer to the Ramond--space of $A$.

\section{Explicit Examples}
\subsection{Self duality of $A_n$}

We consider the example of the Jacobian Frobenius Algebra associated to the function $z^{n+1}$

$$A_n:= {\bf C}[z]/(z^{n})$$

together with the induced multiplication and the Grothendieck residue.
Explicitly:
$$
z^i z^j= \left\{ \begin{tabular} {ll}$z^{i+j}$& \mbox{if $i+j \leq n$}\\
$0$&\mbox{else}\\
\end{tabular}
\right.
$$
the form
$$
\eta(z^i,z^j) = \delta_{i,n-1-j}
$$
and the grading:

$$
E(z^i):= \frac{i}{n+1}
$$
which means $\rho=z^{n-1}$ and $d= \frac{n-1}{n+1}$.

We consider just the group $J \simeq \Znn$ with the generator $j$ acting on $z$ by multiplication
with $\zeta_{n+1}:= \exp(2\pi i\frac{1}{n+1})$. We have
$$Fix_{j^i} \left\{ \begin{tabular} {ll}${\bf C}$& \mbox{if $i=0$}\\
$0$&\mbox{else}\\
\end{tabular}
\right.
$$

and thus
$$
A_{j^i} \left\{ \begin{tabular} {ll}$A_n$& \mbox{if $i=0$}\\
$1_{j^i}$&\mbox{else}\\
\end{tabular}
\right.
$$

Furthermore we have the following grading;

$$
(E, \bar E)(1_{j^i})=
\left\{ \begin{tabular} {ll}
$(0,0)$& for $i=0$\\
$(\frac{k-1}{n+1}, \frac{n-k}{n+1})$&else
\end{tabular}
\right.
$$

which means $\rho_{j^i}=1_{j^i}$ and $d_{j^i}= 0$.

Using the reconstruction Theorem we have to find a cocycle $\g$ and a compatible action $\varphi$.
There is no problem for the metric since always $|N_g|=1$ and if
$n+1$ is even $\det(\zeta^{\frac{n+1}{2}})=-1$.
Since the group $J$ is cyclic there is no freedom of choice for $\eps$ and just two choices of
parity are possible corresponding to $j\mapsto \pm 1$.

From the general considerations we know $\g_{j^i,j^{n-1-i}}\in A_e$
and $\deg(\g_{j^i,j^{n-1-i}}) = d-d_{j^i}= \frac{n-1}{n+1}$ which yields
$$
\g_{j^i,j^{n-1-i}}=( (-1)^{\tilde j i}\zeta^{i})^{1/2}\rho = ( (-1)^{\tilde j}\zeta)^{i/2} z^{n-1}
$$
for the other $\gamma$ notice that $\deg(1_{j^i}) +\deg (1_{j^k})= \frac{i+k-2}{n+1}$ while
$\deg(1_{j^{i+j}})= \frac{i+k-1}{n+1}$ if $i+k\neq n+1$, but there is no element
of degree $\frac{1}{n+1}$ in $A_{j^{i+k}}$ for $ i+k\neq n+1$.

Hence
$$
\gamma_{j^i,j^k}=
\left\{ \begin{tabular} {ll}
$((-1)^{\tilde j}\zeta)^{i/2}z^{n-1}$& for $i+j=n+1$\\
$0$&else
\end{tabular}
\right.
$$
this means that
$$
\varphi_{j^i,j^k}= (-1)^{\tilde j^i \tilde j^k}\zeta^{-i}
$$

Therefore the $G$--invariants $A^G=A_1$  are generated by the identity $1$.

The Ramond grading of this algebra is
$$
(E, \bar E)(1_{j^i}v)=
\left\{ \begin{tabular} {ll}
$(0,0)$& for $i=0$\\
$(\frac{k}{n+1}-\frac{1}{2}, -\frac{k}{n+1}+\frac{1}{2})$&else
\end{tabular}
\right.
$$

Since $j\in G$ the special $G$--Frobenius algebra is Euler and the dual is defined, moreover $G=J$
so that the vector--space structures of $A$ and $\check A$ coincide.
The grading is given by:
$$
(E, \bar E)(1_{j^i}v)=
\left\{ \begin{tabular} {ll}
$(0,0)$& for $i=0$\\
$(-\frac{k-(n+1)}{n+1}-\frac{1}{2}, \frac{n+1-k}{n+1}+\frac{1}{2})$&else
\end{tabular}
\right.
$$

The $G$--action is given by
$$
\check\varphi_{j^i,j^k}=
\left\{
\begin{tabular} {ll}
$(-1)^{\tilde j^i \tilde j^k}\zeta^{i}$& for $k=1$\\
$(-1)^{\tilde j^i \tilde j^k}$& else
\end{tabular}
\right.
$$
This $G$--action leaves all even sectors $\check A_{j^i}$ invariant except for $i=1$
if $\tilde j =1$.
Thus with the choice of all even sectors
we have as $G$--modules

$\check A \simeq A_n$

where more explicitly

$\check 1_{j^i} \mapsto z^{n+1-k}:k=2,\dots n$
and $\check 1_0 \mapsto 1$.

Notice that this $A_n$ is of $(a,c)$--type however.
Also since $J=G$ the form $\bar\eta$ pulls back and gives a non--degenerate form on $\check A^G$,
which is the usual form on $A_n$.
Furthermore the usual multiplication on $A_n$ is compatible with everything so that
we can say that $A_n$ is self--dual under this operation.

\subsection{$D_n$ from a special $\Z2$--Frobenius algebra based on $A_{2n-3}$}
In this section, we show how to get $D_n$ from a special $\Z2$--Frobenius algebra based on $A_{2n-3}$.
The function for the Frobenius algebra $A_{2n-3}$ is $z^{2n-2}$. Since this is an even
function,
$\Z2$ acting via $z\mapsto -z$ is a symmetry. There are two sectors, the untwisted and the
twisted sector containing the element $1_{-1}$ with degree $0$.
The multiplication is fixed by $\deg(\g_{-1,-1})=\frac{2n-4}{2n-2}$ thus
$$
\g_{-1,-1}= z^{2(n-2)}
$$
again the group is cyclic so the $G$--action only depends on the choice of parity of the $-1$--sector.

In the untwisted sector we have $A_e^{\Z2}= \la 1,z^2,\dots,z^{2(n-1)}\ra \simeq A_{n-1}$
and the action of $-1$ on $1_{-1}$ is given by
$$
\varphi_{-1,-1}= (-1)^{\widetilde{-1} \widetilde {-1} +1}
$$
so that if  $\sigma(-1) =1$
$$
A_{2n-3,\Z2}^{\Z2} \simeq
D_n
$$
and if $\sigma(-1)=0$ we just obtain the invariants of the untwisted sector
which are isomorphic to $A_{n-1}$.

The untwisted sector is given by the singularity $A_{n-1}$ as expected upon the transformation
$u= z^{2}$.
Notice that the invariants of the Ramond sector yield the singularity $A_{n-2}$
as expected from [Wa].
These are of the form $u^i du$ or $z^{2i+1}dz$ with $i=0,\dots, n-3$.

\subsection{Point mod $G$}
In the theory of Jacobian Frobenius algebras there is the notion of a point played by a Morse
singularity $z_1^2 +\dots +z_n^2$.
Any finite subgroup $G \subset O(n,k)$ leaves this point invariant.

The $G$--twisted algebra after possibly stabilizing is the following.
$$
A= \bigoplus_{g\in G} k 1_g
$$
And the grading is $\deg(1_g) = (\frac{1}{2}s^-_g,\frac{1}{2}s^-_g)$, since $d=d_g=1$.

Using \ref{zerocheck} it follows  that the $\psi$ cannot vanish,
thus fixing $\phi$ and $\eps$, so that the possibilities are
enumerated by the graded cocycles. The compatibility equations hold 
automatically.

{\sc Special case.}

If we assume that $G\subset O(n,{\bf C})$ and that $s^-_g=0$
(i.e.\ $\sum_{i:\l_i \neq 0} \frac{1}{2\pi i} \l_i
=\frac{|N_g|}{2} \in {\bf N}$), then we the cocycles lie in
$H^2(G,k^*)$ and the possible algebra structures are those of
twisted group algebras.

\subsubsection{Point mod $\Zn$}

By the above analysis we realize $\Zn$ as the sub--group of rotations of order $n$ in ${\bf C}$.
We have that $s_g^- = 0$ and thus we can choose the full cocycle making
$A$ into $A_n$, multiplicatively, with trivial
grading and trivial $G$--action if one chooses all even sectors.
The metric, however, will not be consistent with $A_n$.
The identity pairs with itself for instance. Dualizing $A$, we obtain the following space
$$
\check A= \oplus_{j^i} A_{j^i-1}
$$
with again a trivial $(E,\bar E)$ grading.
Choosing the generator $\tilde j := j^{-1}$ for $J$, the metric reads
$$\tilde\eta(\tilde 1_{\tilde j^i},\tilde 1_{\tilde j^k})=
\left\{
\begin{array}{ll}
1 & \mbox{if $ i+j =n-1$}\\
1 & i=j =n\\
0&\mbox{else}
\end{array}
\right.
$$
This metric is compatible with the following multiplication:
$$
\tilde 1_{\tilde j^i}\circ \tilde 1_{\tilde j^k}=
\left\{
\begin{array}{ll}
 \tilde 1_{\tilde j^{i+k}} & \mbox{if $ i+j \leq n-1$}\\
 \tilde 1_{\tilde j^{n-1}} & i=j =n\\
0&\mbox{else}
\end{array}
\right.
$$

The structure of this algebra is the one obtained in the $r$--spin
curve picture of the $A_{n}$--model. Which in mirror symmetry
parlance the A--model version of the theory. By our mirror symmetry we
see that we obtain this A--model by orbifolding the Landau--Ginzburg
B--model pair $(A_{n},A_{1})$ above by $\Zn$. We previous did this by
regarding the left side of the pair and the above calculation is for

the right side of the pair. Thus we see that the appearance of the
extra basis element can be seen as natural from this point of view.
Answering the question of its emergence in [W].

\end{document}